\documentclass[11pt,twoside]{atmp}

\usepackage{amsmath,amssymb}

\usepackage[all]{xy}
\usepackage{color}

\newcommand{\ZZ}{{\bf Z}}

\newcommand{\RR}{{\bf R}}
\newcommand{\CC}{{\bf C}}

\newcommand{\FF}{{\bf F}}
\newcommand{\HH}{{\bf H}}
\newcommand{\PP}{{\bf P}}

\makeatletter
\renewcommand{\@begintheorem}[2]{
\rm \trivlist \item [\hskip \labelsep {\bf #2\ \ #1.}]
                                }
\makeatother

\makeatletter

\DeclareFontFamily{U}{cyr}{}
\DeclareFontShape{U}{cyr}{m}{n}{
  <5> wncyr5 <6> wncyr6 <7> wncyr7 <8> wncyr8 <9> wncyr9 <10->
wncyr10}{}
\DeclareMathAlphabet{\mathcyr}{U}{cyr}{m}{n}

\begin{document}

\title[Siegel modular forms and finite symplectic groups]
{Siegel modular forms and finite symplectic groups}


\author[F. Dalla Piazza, B. van Geemen]
{Francesco Dalla~Piazza$^1$, Bert van Geemen$^2$}

\address{$^1$Dipartimento di Scienze Fisiche e Matematiche, \\
Universit\`a dell'Insubria,\\
Via Valleggio 11, I-22100 Como, Italia. \\ [0.3em]
$^2$Dipartimento di Matematica,\\
Universit\`a di Milano,\\
Via Saldini 50, I-20133 Milano, Italia.}
 
\addressemail{$^1$f.dallapiazza@uninsubria.it,
$^2$geemen@mat.unimi.it}

\begin{abstract}
The finite symplectic group $Sp(2g)$ over the field of two elements 
has a natural representation on the vector space of Siegel modular forms of given weight for the principal congruence subgroup of level two. In this paper we decompose this representation, for various (small) values of the genus and the level, into irreducible representations. As a consequence we obtain uniqueness results for certain modular forms related to the superstring measure, a better understanding of certain modular forms in genus three studied by D'Hoker and Phong as well as a new construction of Miyawaki's cusp form of weight twelve in genus three. 
\end{abstract}

\maketitle

\section*{Introduction}

The vector spaces of Siegel modular forms for the principal
congruence subgroup $\Gamma_g(n)$ of level $n$ are of great interest in various areas in mathematics. For example in geometry, where they correspond to holomorphic tensors on moduli spaces of abelian varieties, in arithmetic, by way of the L-series of Hecke eigenforms and their relation to Galois representations, and in mathematical physics, where they come up in the bosonic and superstring measures. As the dimension of these spaces grows fast with both genus $g$ and level $n$ it is convenient to use the natural
representation of the finite symplectic group $\Gamma_g/\Gamma_g(n)$
on these spaces. This allows one to pick out, in a canonical fashion,
subspaces, i.e.\ subrepresentations, which can then be analyzed more thoroughly. 

In this paper we will be concerned only with the subgroup $\Gamma_g(2)$
and its quotient $Sp(2g)$, the  finite symplectic group 
over the field of two elements. For genus $g\leq 3$ the spaces of 
modular forms on $\Gamma_g(2)$ of even weight have an explicit description in terms of theta series and this allows us to determine 
corresponding representations explicitly for small $k$, 
usually with the help of a computer. In general, the theta series generate only a subspace of the space of all modular forms, but as this subspace is an $Sp(2g)$-subrepresentation, it is still interesting to find its decomposition. 

As an application of our results we prove in section \ref{DP}
some uniqueness results on the modular forms $\Xi_8[0^{(g)}]$ for small $g$. These modular forms 
are of considerable interest in superstring theory 
(see \cite{Mo}). We also investigate the space spanned by the $336$ modular forms of weight $6$ and genus three
studied extensively by D'Hoker and Phong in \cite{DP2}. We show that this space has dimension $105$, that the 
$Sp(6)$-representation on it is irreducible and we 
identify the representation in the character table of Frame 
\cite{Frame}. Moreover, any function in this space turns out to be a cusp form and the sum of the squares of the 336 functions is the cusp form $F_{12}$ found by Miyawaki whose Spinor L-function was recently  determined by B.\ Heim.

In appendix B (section \ref{B}) we recall the basics on the transformation theory of the theta functions which we used in order to determine the representations of $Sp(2g)$ on modular forms.



\noindent

\section{Siegel modular forms of level two}
\label{Siegel modular forms}

\subsection{} In this section we recall the action of a finite symplectic group 
on the modular forms of a given weight $k$ for the congruence subgroup of $\Gamma_g(n)$ of $Sp(2g,\ZZ)$ in the case $n=2$ (for other positive integers  
one has a similar action).

\subsection{The groups $\Gamma_g$, $\Gamma_g(n)$ and $Sp(2g)$.}
We will write $\Gamma_g:=Sp(2g,\ZZ)$ and
$$
\Gamma_g(n)\,:=\,\{\,M\in Sp(2g,\ZZ):\; M\,\equiv\, I\;\mbox{mod}\,n\,\},
$$
is a normal subgroup. In case $n=2$, we have 
$$
\Gamma_g(2)=\ker(Sp(2g,\ZZ)\,\longrightarrow\,Sp(2g)\,:=\,Sp(2g,\FF_2))
$$
where we write $\FF_2:=\ZZ/2\ZZ$ for the field with two elements.
The reduction mod two map is surjective 
(\cite{Igusa}, V.6 Lemma 25)
and thus we have an isomorphism
$$
Sp(2g)\,\cong\,\Gamma_g/\Gamma_g(2).
$$
The order of this group is (cf.\ \cite{Jac}, Thm 6.18):
$$
|Sp(2g)|\,=\,2^{2g-1}(2^{2g}-1)|Sp(2g-2)|\,=\,2^{2g-1}(2^{2g}-1)2^{2g-3}(2^{2g-2}-1)
\cdots 2(2^2-1).
$$
In particular, $|Sp(2g)|=6,6!=720,36\cdot (8!)$ for $g=1,2,3$.

The group $Sp(2g)$ acts (linearly) on the $2^{2g}$ points of $\FF_2^{2g}$ (sometimes called period characteristics) and it acts non-linearly on the (theta) characteristics $[{}^{a_1\ldots a_g}_{b_1\ldots b_g}]$ with $a_i,b_i\in\{0,1\}$, see Appendix A. 
It has two orbits on the characteristics, of length $2^{g-1}(2^g+1)$ 
and $2^{g-1}(2^g-1)$ respectively, the first orbit consists of the even characteristics which are those with $\sum a_ib_i\equiv 0$ mod $2$,
the other orbit are the odd characteristics.

One has the following description of $Sp(2g)$ for small $g$:
the action of $Sp(2)$ on the three non-zero points in $\FF_2^2-\{0\}$
induces an isomorphism
$Sp(2)=SL(2,\FF_2)\cong S_3$, 
the action of $\Gamma_4$ on the six odd characteristics (for $g=2$)
induces an isomorphism
$Sp(4)\cong S_6$,
and $Sp(6)$ is isomorphic to the subgroup of elements with determinant $1$ of
the Weyl group $W(E_7)$ in its standard 7-dimensional representation
(\ref{transvections}, \cite{Atlas} where $Sp(6)=S_6(2)$).

\subsection{The subgroups $O^+(2g)$ and $O^{-}(2g)$ of $Sp(2g)$}\label{O+,O-}\label{Gamma(2n,4n)}
The group $\Gamma_g$ acts on $V=\FF_2^{2g}$, and thus on the characteristics,
through its quotient $Sp(2g)=\Gamma_g/\Gamma_g(2)$. The stabiliser subgroup of the
even characteristic $[0]:=[{}^{0\ldots0}_{0\ldots0}]$ is the subgroup $\Gamma_g(1,2)$, where we define more generally subgroups
(cf.\ \cite{Igusa}, V.2, p.178):
$$
\Gamma_g(n,2n)\,:=\,
\{M\in \Gamma_g(n):\,
{\rm diag}A\,{}^t\!B \equiv {\rm diag }C\,{}^t\!D\equiv 0\,{\rm mod}\,2n\,\}.
$$
In case $n$ is even, $\Gamma_g(n,2n)$ is a normal subgroup of $\Gamma_g(1)$.

The image of $\Gamma_g(1,2)$ in $Sp(2g)$ is denoted by $O^+(2g)$:
$$
O^+(2g)\,:=\,\Gamma_g(1,2)/\Gamma_g(2)\qquad(\subset Sp(2g)).
$$
As $Sp(2g)$ acts transitively on the even theta characteristics, there is a natural bijection 
$$
Sp(2g)/O^+(2g)\,\longrightarrow\,\{\Delta:\,\Delta\;\mbox{even}\,\},\qquad
hO^+(2g)\,\longmapsto\, h\cdot [0].
$$ 
In particular, $[Sp(2g):O^+(2g)]=2^{g-1}(2^g+1)$. 
One has $O^+(2)\cong\ZZ/2\ZZ$ and $O^+(4)$ is isomorphic to the subgroup of the symmetric group $S_6\cong Sp(4)$ of permutations 
$\sigma$ such that $\sigma(\{1,2,3\})\subset \{1,2,3\}$ or
$\sigma(\{1,2,3\})= \{4,5,6\}$, thus $S_3\times S_3$ is a subgroup of index two in $O^+(4)$ and
$|O^+(4)|=(3!)\cdot(3!)\cdot2=72$.
One has $O^+(6)\cong S_8$,
the symmetric group of order $8!$ and $O^+(8)$ is the quotient of the subgroup of elements of $W(E_8)$ with determinant $+1$ in the standard $8$-dimensional representation, by the center, generated by $-I$. 

Similarly, we define a subgroup $O^-(2g)\subset Sp(2g)$ as the stabilizer of
the odd characteristic $[{}^{10\ldots0}_{10\ldots0}]$ and we have $[Sp(2g):O^+(2g)]=2^{g-1}(2^g-1)$. One has $O^-(2)\cong Sp(2)\cong S_3$,
$O^-(4)\cong S_5$, a symmetric group, and $O^-(6)\cong W(E_6)$, the Weyl group of the root system $E_6$.

\subsection{Modular forms of level $2$}\label{level 2}
A Siegel modular form of genus $g$, weight $k$ and level $2$ is a holomorphic function on the Siegel upper half space of genus $g$
$$
f:\HH_g\,\longrightarrow\,\CC,\qquad
f(M\cdot\tau)=\det(C\tau+D)^kf(\tau)\qquad \forall M\in \Gamma_g(2)
$$
(and in case $g=1$ one should also impose a growth condition on $f$).
Here, as usual, the action of $\Gamma_g$ on the Siegel upper half space 
is given by
$$
M\cdot \tau\,:=\,(A\tau+B)(C\tau+D)^{-1}\qquad
M:=\begin{pmatrix}A&B\\C&D\end{pmatrix}\in Sp(2g,\ZZ),\quad\tau\in\HH_g.
$$

The Siegel modular forms of genus $g$, weight $k$ and level $2$
form a finite dimensional complex vector space 
denoted by $M_k(\Gamma_g(2))$.
The finite group $Sp(2g)$ has a representation 
$$
\rho=\rho_k:Sp(2g)\,\longrightarrow \,GL(M_k(\Gamma_g(2)))
$$ 
on this vector space defined by
$$
(\rho(g^{-1})f)(\tau)\,:=\,\det(C\tau+D)^{-k}f(M\cdot \tau),
$$
where
$M\in \Gamma_g$ is a representative of $g\in Sp(2g)$ and $f\in M_k(\Gamma_g(2))$
(note that 
$\det(C\tau+D)^{-k}f(M\cdot \tau)=f(\tau)$ for $M\in \Gamma_g(2)$,
thus the action of $M\in \Gamma_g$ factors over $\Gamma_g/\Gamma_g(2)=Sp(2g)$).
The equality $\rho(gh)=\rho(g)\rho(h)$ for $g,h\in Sp(2g)$ follows from 
$(MN)\cdot \tau=M\cdot(N\cdot \tau)$  and
$\gamma(MN,\tau)=\gamma(M,N\cdot \tau)\gamma(N,\tau)$
where $\gamma(M,\tau):=\det(C\tau+D)$.

\subsection{The theta constants $\Theta[\sigma]$} \label{Theta[sigma]}
To determine the modular forms of even weight on $\Gamma_g(2)$ it is convenient to
define the $2^g$ theta constants (see Appendix A, \ref{thetas}): 
$$
\Theta[\sigma](\tau)\,:=\, \theta[{}^\sigma_{0}](2\tau,0),\qquad
[\sigma]=[\sigma_1\;\sigma_2\;\ldots\;\sigma_g],\;\sigma_i\in\{0,1\},\;
\tau\in\HH_g.
$$
These theta constants are "modular forms of weight $1/2$"  on the normal subgroup 
$\Gamma_g(2,4)\subset \Gamma_g(2)$, defined in \ref{O+,O-}.

Basically, the invariants of degree $4k$ of the quotient group $\Gamma_g(2)/\Gamma_g(2,4)\cong\FF_2^{2g}$ in the ring of polynomials in the $\Theta[\sigma]$'s are modular forms of weight $2k$ on $\Gamma_g(2)$. 
Due to the half integral weight, the quotient group $\Gamma_g(2)/\Gamma_g(2,4)$ doesn't actually act on the $\Theta[\sigma]$, 
but a central extension, the Heisenberg group $H_g$ does (see Appendix B) and we have to take the invariants for this Heisenberg group.
The subspace of $M_{2k}(\Gamma_g(2))$ of these Heisenberg invariants
is denoted by: 
$$
M_{2k}^\theta(\Gamma_g(2))\,\subset\,M_{2k}(\Gamma_g(2)),\qquad
M_{2k}^\theta(\Gamma_g(2))\,:=\,\CC[\ldots,\Theta[\sigma],\ldots]_{2k}^{H_g}
$$
where $\CC[\ldots,\Theta[\sigma],\ldots]_{2k}$ is the subspace of homogeneous polynomials of degree $2k$ in the $\Theta[\sigma]$'s.

\subsection{The Heisenberg group $H_g$}\label{Heisenberg}
We define a finite (Heisenberg) group by:
$$
H_g=\mu_4\times \FF_2^g\times\FF_2^g,\qquad
(s,x,u)(t,y,v)=(st(-1)^{uy},x+y,u+v)
$$
where $\mu_4=\{z\in\CC:z^4=1\}$ is the multiplicative group of fourth roots of unity, 
and $uy=\sum u_iy_i$ for $u=(u_1,\ldots,u_g), y=(y_1,\ldots,y_g)\in\FF_2^g$.
The center of $H_g$ is $\mu_4$ and the quotient $H_g/\mu_4$ is isomorphic to $\FF_2^{2g}$. 

Next we define the (Schr\"odinger) representation of $H_g$ on polynomials in the
$2^g$ variables $X_\sigma$, where $[\sigma]=[\sigma_1\;\sigma_2\;\ldots\;\sigma_g],\;\sigma_i\in\{0,1\}$:
$$
(s,x,u)X_\sigma\,:=\,s(-1)^{(x+\sigma)u}X_{x+\sigma},
$$
where we now consider $\sigma_i\in\FF_2=\{0,1\}$, so $x+\sigma=[x_1+\sigma_1,\ldots]$.
This action is extended in the obvious way to polynomials in the $X_{\sigma}$'s.
The action of $H_g$ on $\CC[\ldots,X_\sigma,\ldots]$ induces the one on
$\CC[\ldots,\Theta[\sigma],\ldots]$ under the map $X_\sigma\mapsto \Theta[\sigma]$,
(c.f.\ Appendix B).

\subsection{Heisenberg invariants}
Note that a homogeneous polynomial in the $X_{\sigma}$ is invariant under the subgroup
$\mu_4$ of $H_g$ iff its degree is a multiple of $4$, so $H_g$-invariant polynomials are modular forms of even weight. 

The subring of Heisenberg invariant polynomials in 
$\CC[\ldots,X_{\sigma},\ldots]$ will be denoted by $\CC[\ldots,X_{\sigma},\ldots]^{H_g}$. 
In section \ref{dimension formula} we give a formula for the dimension of the vector space of Heisenberg invariant polynomials of degree $4n$.

It is not hard to construct a basis of this vector space.
Each monomial in a Heisenberg invariant polynomial $P$ is a product $\prod_{i=1}^{4n}X_{\sigma_i}$ which must be invariant under
the action of the $(1,0,u)\in H_g$, hence $\sum \sigma_i=0$.
As $P$ is also invariant under the $(1,x,0)$, 
all monomials of the type $\prod_{i=1}^{4n}X_{\sigma_i+x}$, for any $x\in \FF_2^g$, have the same coefficient in $P$. 
Thus a basis is obtained by the $\sum_x\prod_{i=1}^{4n}X_{\sigma_i+x}$, where $\sum \sigma_i=0$.

\subsection{A dimension formula}\label{dimension formula}
We give a formula for the dimension of the Heisenberg invariants
in the vector space 
$
\CC[\ldots,X_{\sigma},\ldots]_n
$
of homogeneous polynomials of degree $n$ in the $X_{\sigma}$'s.
Let, for fixed $g$,
$$
\rho_n\,:\,H_g\longrightarrow\, GL(\CC[\ldots,X_{\sigma},\ldots]_n)
$$
be the representation of $H_g$ on this vector space. 
Then the space of Heisenberg invariants is the trivial subrepresentation of $\rho_n$, so its dimension is
$$
\dim \,\CC[\ldots,X_{\sigma},\ldots]_n^{H_g}\,=\,
\langle \rho_n,1_{H_g}\rangle_{H_g},
$$
the multiplicity of the trivial representation $1_{H_g}$ of $H_g$ in $\rho_n$. This scalar product of characters is given by
$$
\langle \rho_n,1_{H_g}\rangle_{H_g}\,=\,
\mbox{$\frac{1}{|H_g|}$}\sum_{h\in H_g} tr(\rho_n(h))
$$
where $tr$ is the trace.

Note that if $xu=0$, 
but $(x,u)\neq (0,0)$, $(1,x,u)$ has order two in $H_g$ and the
eigenvalues of $(t,x,u)\in H_g$ on $\CC[\ldots,X_{\sigma},\ldots]_1$
are $t$ and $-t$, each with multiplicity $2^{g-1}$ for all $t\in\mu_4$. 
In case $xu=1$, $(1,x,u)$ has order four and $(t,x,u)$ has eigenvalues
$it$ and $-it$ (with $i^2=-1$), each with multiplicity $2^{g-1}$.
So if $(x,u)\neq (0,0)$ the eigenvalues of $\rho_1(t,x,u)$, with $t\in\mu_4$, are 
$i^a,-i^a$ with $a=1$ or $a=2$, each with multiplicity $2^{g-1}$. 

If $\alpha_1,\ldots,\alpha_N$, $N=2^g$ are the eigenvalues of 
$(t,x,u)$ on $\CC[\ldots,X_{\sigma},\ldots]_1$, the eigenvalues
on $\CC[\ldots,X_{\sigma},\ldots]_n$ are the $\alpha_1^{m_1}\alpha_2^{m_2}\ldots\alpha_N^{m_N}$ with $\sum m_i=n$.
As the trace is the sum of the eigenvalues, we get (with a variable $X$)
$$
\sum_n tr(\rho_n(t,x,u))X^n\,=\,\prod_{i=1}^N (1-\alpha_iX)^{-1}.
$$
So if $(x,u)\neq (0,0)$ we have 
$$
\sum_n tr(\rho_n(t,x,u))X^n\,=\,(1-i^{2a}X^2)^{-2^{g-1}}
\,=\,\sum_m (-1)^{am}\binom{2^{g-1}+m-1}{m}X^{2m},
$$
and in case $(x,u)=(0,0)$ the trace is just
$$
\sum_n tr(\rho_n(t,0,0))X^n\,=\,
\sum_nt^n(\dim\CC[\,\ldots,X_{\sigma},\ldots]_n)X^n
\,=\,\sum_n t^n\binom{2^{g}+n-1}{n}X^n.
$$
This leads to the formula (note that there are non-trivial invariants only if $n$ is a multiple of $4$):
$$
\dim (\CC[\,\ldots,X_{\sigma},\ldots]_{4n})^{H_g}\,=\,
2^{-2g}\left(\binom{2^g+4n-1}{4n}+
(2^{2g}-1)\binom{2^{g-1}+2n-1}{2n}\right).
$$
For small $g$ we list some of these dimensions in the table below.
$$
\begin{array}{|@{\hspace{6pt}}l||@{\hspace{6pt}}c|@{\hspace{6pt}}c|
@{\hspace{6pt}}c|@{\hspace{6pt}}c|}
\hline g\,\backslash\,\mbox{degree}& 
4&
8&
12&
16
\\\hline \hline
1&2&3&4 &5\\ \hline
2&5&15&35 &69\\ \hline
3&15&135&870&3993\\ \hline
4&51&2244&69615&1180395\\ \hline
\end{array}
$$
In case $n=1$, the dimension formula gives:
$$
\dim (\CC[\,\ldots,X_{\sigma},\ldots]_{4})^{H_g}\,=\,
(2^g+1)(2^{g-1}+1)/3.
$$

\

\subsection{The ring of modular forms}\label{ring of modular forms}
The vector spaces $M_{2k}^\theta(\Gamma_g(2))$ are images of the spaces of $H_g$ invariants of degree $4k$ under the surjective maps
$$
\CC[\ldots,X_\sigma,\ldots]_{4k}^{H_g}\,\longrightarrow \, M_{2k}^\theta(\Gamma_g(2)):=\CC[\ldots,\Theta[\sigma],\ldots]^{H_g}_{4k},
\qquad X_\sigma\,\longmapsto\,\Theta[\sigma].
$$
These maps define a surjective $\CC$-algebra homomorphism 
$$
\CC[\ldots,X_\sigma,\ldots]^{H_g}\,\longrightarrow\,
M^\theta(\Gamma_g(2))\,:=\,\oplus_k M_{2k}^\theta(\Gamma_g(2))
$$
whose kernel is the ideal of algebraic relations between the $\Theta[\sigma]$'s.
So a polynomial $F(\ldots,X_\sigma,\ldots)$ maps to zero iff $F(\ldots,\Theta[\sigma](\tau),\ldots)=0$ for all $\tau\in\HH_g$.
In geometrical terms, a homogeneous polynomial $F$ lies in the kernel iff
the image of $\HH_g$ under the map
$$
\HH_g\,\longrightarrow\,\PP^{2^g-1},\qquad
\tau\,\longmapsto\,(\ldots:\Theta[\sigma](\tau):\ldots)
$$
lies in the zero locus $Z(F)\subset\PP^{2^g-1}$. We will describe this map for $g=1,2,3$, in particular the image of the map contains an open subset for $g=1,2$, thus there are no polynomials vanishing on the image. 
In case $g=3$ the image is a (Zariski) open subset of a hypersurface $Z(F_{16})\subset\PP^7$, 
for a certain homogeneous polynomial $F_{16}$ in eight variables, 
hence $M^\theta(\Gamma_g(2))\cong \CC[\ldots,X_\sigma,\ldots]^{H_g}/(F_{16})$.

For any $g$, $g\geq 4$, there must be (many) algebraic relations between the $\Theta[\sigma]$'s because $\dim \HH_g = g(g+1)/2 < 2^{g}-1$, but a complete description of these relations is not known.

The graded ring of modular forms of even weight on $\Gamma_g(2)$ is the normalization of the ring of
the $\Theta[\sigma]$'s (cf.\ \cite{SM2} Thm 2, \cite{R1}, \cite{R2}):
$$
\oplus_{k=0}^\infty\,M_{2k}(\Gamma_g(2))\;=\;
(\CC[\ldots,\Theta[\sigma],\ldots]^{H_g})^{Nor}.
$$

In case $g=1,2$ there are no relations and the ring of invariants is already normal. In case $g=3$, there is one relation $F_{16}(\ldots,\Theta[\sigma],\ldots)=0$ homogeneous of degree 16, and Runge (\cite{R1}, \cite{R2}) showed that the quotient of the ring of invariants by the ideal generated by this relation is again normal. 
This implies that any modular form of weight $2k$ can be written as a homogeneous polynomial of degree $4k$ in the
$\Theta[\sigma]$'s if $g\leq 3$, thus
$$
M_{2k}^\theta(\Gamma_g(2))\,=\,M_{2k}(\Gamma_g(2))\qquad
\mbox{for}\quad g=1,2,3.
$$

This polynomial is unique for $g\leq 2$.
For $g=3$ it is unique if its degree is less then 15 and else it is
unique up to the addition of $F_{16}G_{4k-16}$ where $G_{4k-16}$
is any homogeneous polynomial of degree $4k-16$ in the $\Theta[\sigma]$'s. 

For $g>3$ there will always be non-trivial relations 
and if $g>4$ the ring 
$\CC[\ldots,\Theta[\sigma],\ldots]^{H_g}$'s is not normal,
(cf.\ \cite{OSM}, Theorem 6, but note that our $H_g$ is different from theirs).
In case the ring is not normal, there are also quotients $G_{4k+d}/H_{d}$
of homogeneous polynomials in the $\Theta[\sigma]$'s, of degree $4k+d$ and $d$ respectively, which are modular forms of weight $4k$ (but which cannot be written as a polynomial in the $\Theta[\sigma]$'s).

\section{Classical theta constants and representations}
\label{small mod}

\subsection{} To describe the spaces of modular forms $M_{2k}^\theta(\Gamma_g(2))$ it is convenient to use also the classical theta functions with (`half integral') characteristics $\theta[\Delta]$.
In particular, we are interested in decomposing these spaces
into irreducible representations for the group $Sp(2g)$ and we want to describe the 
subspace of $O^+$-invariants as well as $O^+$-anti-invariants (cf.\ \ref{O+ inv}),
these will have applications to the superstring measures.

\subsection{The quadratic relations between the $\theta[\Delta]$'s and the $\Theta[\sigma]$'s} A classical formula for theta functions shows that any product of two $\Theta[\sigma]$'s is a linear combination of the $\theta[\Delta]^2$.
Note that there are $2^g$ functions $\Theta[\sigma]$ and thus there are
$(2^g+1)2^g/2=2^{g-1}(2^g+1)$ products $\Theta[\sigma]\Theta[\sigma']$.
This is also the number of even characteristics, and the products
$\Theta[\sigma]\Theta[\sigma']$ span the same space (of modular forms of weight 1) as the $\theta[\Delta]^2$'s, which has dimension $2^{g-1}(2^g+1)$.

As the degree of an $H_g$-invariant homogeneous polynomial in the $\Theta[\sigma]$ is a multiple of four, say $4k$, it can be written as a homogeneous polynomial of degree $2k$ in the $\theta[\Delta]^2$'s. Thus for $g\leq 3$, any element in
$M_{2k}(\Gamma_g(2))$ is a homogeneous polynomial of degree $2k$ in the
$\theta[\Delta]^2$'s. 

The $\theta[\Delta]^2$ are the better known functions and their transformation under $\Gamma_g(1)$ is easy to understand, but
the $\Theta[\sigma]$ have the advantage that they are algebraically independent for $g\leq 2$ and there is a unique relation of degree $16$ for $g=3$. In contrast, there are many quadratic relations
between the $\theta[\Delta]^2$'s, for example Jacobi's relation in $g=1$.

\subsection{A classical formula}\label{classical}
The classical formula used here is (cf.\ \cite{Igusa} IV.1):
$$
\theta[{}^a_{b}]^2\,=\,\sum_{\sigma}
(-1)^{\sigma b}\Theta[\sigma]\Theta[\sigma+a]
$$
where we sum over the $2^g$ vectors $\sigma\in\FF_2^g$ and $[{}^a_{b}]$
is an even characteristic, so $a\,{}^t\!b=0\;(\in\FF_2)$.
These formulae are easily inverted to give:
$$
\Theta[\sigma]\Theta[\sigma+a]\,=\,\mbox{$\frac{1}{2^g}$}\sum_{b}
(-1)^{\sigma b}\theta[{}^a_{b}]^2.
$$
It is easy to see that the $\theta[\Delta]^2$ span one-dimensional subrepresentations of $H_g$. In fact, using the classical formula one finds:
$$
(s,x,u)\theta[{}^a_{b}]^2=
s^2(-1)^{ua +xb}\theta[{}^a_{b}]^2.
$$
This implies that the $\theta[{}^a_{b}]^4$ are 
Heisenberg invariants and thus are in $M_2(\Gamma_g(2))$.
More generally, we have: 
$$
\prod_i^{2k}\theta[{}^{a_i}_{b_i}]^2\,\in\,
M_{2k}(\Gamma_g(2))\quad
\mbox{iff}\quad
\sum a_i=\sum b_i=0\;(\in\FF_2).
$$

For example in case $g=1$ one has
$$
\theta[{}^0_0]^2=\Theta[0]^2+\Theta[1]^2,\qquad
\theta[{}^0_1]^2=\Theta[0]^2-\Theta[1]^2,\qquad
\theta[{}^1_0]^2=2\Theta[0]\Theta[1],
$$
or, equivalently,
$$
\Theta[0]^2=(\theta[{}^0_0]^2+\theta[{}^0_1]^2)/2,\qquad
\Theta[1]^2=(\theta[{}^0_0]^2-\theta[{}^0_1]^2)/2,\qquad
\Theta[0]\Theta[1]=\theta[{}^1_0]^2/2.
$$
Note that upon substituting the first three relations in Jacobi's relation
$\theta[{}^0_0]^4=\theta[{}^0_1]^4+\theta[{}^1_0]^4$ one obtains a trivial 
identity in the $\Theta[\sigma]$'s.

\subsection{The $O^+$-invariants and $O^+$-anti-invariants.} \label{O+ inv}
The function $\theta[0]^4=(\sum_\sigma\Theta[\sigma]^2)^2$ is Heisenberg invariant and thus defines a modular form of weight $2$ on $\Gamma_g(2)$. 
For $g\in O^+(2g)$ we have $g\cdot [0]=[0]$ and the explicit transformation formula
for theta constants shows that $\theta[0]^4$ transforms by a non-trivial character which we denote by $\epsilon$:
$$
\rho(g)\theta[0]^4\,=\,\epsilon(g)\theta[0]^4,\qquad 
\epsilon:\,O^+(2g)\longrightarrow \{\pm 1\}.
$$
For $g\geq 3$, this homomorphism is the only non-trivial one dimensional representation of $O^+(2g)$ and its kernel is a simple group.
In case $g=2$, Thomae's formula for $\theta[0]^4$ implies that $\epsilon$ is the product of the sign character on the subgroup $S_3\times S_3$ of $O^+(4)$  and $\epsilon(g)=1$ if $g=(14)(25)(36)$ where we identify $O^+(4)$ with a subgroup of $S_6$ as in \ref{O+,O-}.

For applications to superstring measures, we will be particularly interested in the subspace of $O^+(2g)$-anti-invariants in weight 6:
$$
M_6(\Gamma_g(2))^{\epsilon}\,:=\,\{f\in M_6(\Gamma_g(2))\,:\;
\rho(g)f=\epsilon(g)f\quad \forall g\in O^+(2g)\,\}
$$
and the space of $O^+$-invariants in weight $8$:
$$
M_8(\Gamma_g(2))^{O^+}\,:=\,\{f\in M_8(\Gamma_g(2))\,:\;
\rho(g)f=f\quad \forall g\in O^+(2g)\,\}.
$$

It should be noted that $Sp(2g)$ permutes the 
$\theta[\Delta]^{4k}\in M_{2k}(\Gamma_g(2))$, up to sign if $k$ is odd.
Thus it is not hard to write down some invariants or anti-invariants, but the problem is to find all of them.

\subsection{The dimensions of the $O^+$-(anti)-invariants}
\label{dim O+ inv}
Once the decomposition of an $Sp(2g)$-representation into irreducibles is known, it is easy to find the dimension of the $O^+$-(anti)-invariants
(although this is a case of overkill, since finding these $O^+$-representations is much easier to do using the action of some generators of these groups). Recently M.\ Oura \cite{Oura}
used the methods from \cite{R1},\cite{R2} 
to determine the dimension of the $O^+$-invariants in $M_{2k}(\Gamma_g(2))$ for small $k$ and $g$.

The dimension of the $O^+$-invariants in $V$ is the multiplicity of the trivial representation $1$ of $O^+$ in the $O^+$-representation $Res^{Sp}_{O^+}(V)$ (the restriction of the representation from $Sp(2g)$ to $O^+(2g)$):
$$
\dim V^{O^+}\,=\,\langle\,Res^{Sp}_{O^+}(V),1\,\rangle_{O^+}\,=\,
\langle\,V,Ind^{Sp}_{O^+}(1)\,\rangle_{Sp}
$$
where the second equality is Frobenius reciprocity.
According to Frame \cite{F-ind} one has:
$$
Ind^{Sp}_{O^+}(1)\,=\,{\bf 1}\,+\,\sigma_\theta,\qquad
\dim \sigma_\theta\,=\,2^{g-1}(2^{g}+1)-1=(2^g-1)(2^g+2)/2,
$$
where ${\bf 1}$ is the trivial representation and $\sigma_\theta$ is an irreducible representation of $Sp(2g)$, note that 
$\dim Ind^{Sp}_{O^+}(1)=[Sp(2g):O^+(2g)]=2^{g-1}(2^g+1)$.
Thus if the multiplicity of ${\bf 1},\sigma_\theta$ in $V$ is $n_1,n_\theta$ respectively, then $\dim V^{O^+}=n_1+n_\theta$.

Similarly, the dimension of the $O^+$-anti-invariants in $V$ is the multiplicity of the representation $\epsilon$ of $O^+$ in the $O^+$-representation $Res^{Sp}_{O^+}(V)$ :
$$
\dim V^{\epsilon}\,=\, \langle Res^{Sp}_{O^+}(V),\epsilon\,\rangle_{O^+}
\,=\,
\langle V,Ind^{Sp}_{O^+}(\epsilon)\,\rangle_{Sp}.
$$
According to Frame \cite{F-ind},
the induced representation has two irreducible components:
$$
Ind^{Sp}_{O^+}(\epsilon)\,=\,\rho_\theta\,\oplus\,\rho_r,\qquad
\left\{\begin{array}{rcl}
\dim\rho_\theta&=&(2^g+1)(2^{g-1}+1)/3,\\
\dim\rho_r&=&(2^g+1)(2^g-1)/3.
\end{array}\right.
$$
Thus if the multiplicity of $\rho_\theta,\rho_r$ in $V$ is $n_\theta,n_r$ respectively, then $\dim V^{\epsilon}=n_\theta+n_r$.

\subsection{The $Sp(2g)$-representation on $M_2^\theta(\Gamma_g(2))$}
\label{weight 2 irr}
The representation $\rho_2$ of $Sp(2g)$ on the subspace $M_2^\theta(\Gamma_g(2)) \subset M_2(\Gamma_g(2))$ was shown to be irreducible and isomorphic to $\rho_\theta$ in \cite{vG}.

We briefly recall the proof. As $\theta[0]^4\in M_2^\theta(\Gamma_g(2))$ and $\theta[0]^4$ is an $O^+(2g)$-anti-invariant, the representation $\epsilon$ of $O^+$ is a summand of the restriction of $\rho_2$ to $O^+$:
$$
1\,\geq \,\langle Res^{Sp}_{O^+}(\rho_2),\epsilon\rangle_{O^+}
\,=\,
\langle \rho_2,Ind^{Sp}_{O^+}(\epsilon)\rangle_{Sp}.
$$
As $Ind^{Sp}_{O^+}(\epsilon)=\rho_\theta\oplus\rho_r$, either $\rho_\theta$ or $\rho_r$ must be a summand of $\rho_2$.
For $g\geq 3$, and using  $\dim\CC[\ldots,X_\sigma,\ldots]_4^{H_g}=(2^g+1)(2^{g-1}+1)/3$ (cf.\ \ref{dimension formula}) we have:
$$
\dim \rho_r
\,>\,  \dim \rho_\theta
\,=\, \dim\CC[\ldots,X_\sigma,\ldots]_4^{H_g} 
\,\geq\, \dim M_2^\theta(\Gamma_g(2)).
$$ 
Therefore
we have the irreducible representation $\rho_2=\rho_\theta$
on $M_2^\theta(\Gamma_g(2))$ and
$\CC[\ldots,X_\sigma,\ldots]_4^{H_g}=M_2^\theta(\Gamma_g(2))$
(so there are no $H_g$-invariant algebraic relations of degree $4$),
moreover $\dim M_2^\theta(\Gamma_g(2))=(2^g+1)(2^{g-1}+1)/3$.

The induced representation $Ind^{Sp}_{O^+}(\epsilon)$ is realized on a vector space of dimension $[Sp(2g):O^+(2g)]=2^{g-1}(2^g+1)$ which has a basis $e_\Delta$ parametrized by the even theta characteristics (or equivalently, by the cosets $Sp(2g)/O^+(2g)$). 
There is an $Sp(2g)$-equivariant map
$$
\Phi\,:\,Ind^{Sp}_{O^+}(\epsilon)\,\longrightarrow\,
M_2^\theta(\Gamma_g(2)),\qquad
e_\Delta=g\cdot e_{[0]}\longmapsto \rho_2(g)\theta[0]^4=\epsilon_\Delta\theta[\Delta]^4,
$$ 
for some  $\epsilon_\Delta\in \{\pm 1\}$.
The subrepresentation
$\rho^\theta$ is mapped onto $M_2^\theta(\Gamma_g(2))$ and the kernel of $\Phi$ is the subrepresentation $\rho_r$.

\subsection{The $Sp(2g)$-representation on $Sym^2(M_2^\theta(\Gamma_g(2)))$}
\label{sym2 dec}
Frame has shown that the $Sp(2g)$-representation $Sym^2(M_2^\theta(\Gamma_g(2)))$
decomposes into irreducible representations as follows:
$$
Sym^2(\rho_\theta)\,=\,{\bf 1}\,+\, \sigma_\theta\,+\,\sigma_c,\qquad
\dim\sigma_c\,=\,2^{g-2}(2^g+1)(2^g-1)(2^g+2),
$$ 
and $\sigma_\theta$ as in \ref{dim O+ inv}.
The functions $\theta[\Delta]^8$ are in $M_4^\theta(\Gamma_g(2))$.
They are permuted (without signs) by $Sp(2g)$ and span
the subrepresentation ${\bf 1}+ \sigma_\theta$.
The trivial subrepresentation in $Sym^2(M_2^\theta(\Gamma_g(2)))$
is then spanned by the invariant $\sum_\Delta\theta[\Delta]^8$.
It is easy to verify that the dimension of the image of 
$Sym^2(M_2^\theta(\Gamma_g(2)))$ is larger than
$2^{g-1}(2^g+1)=\dim ({\bf 1}+ \sigma_\theta)$ for $g\geq 2$.
As the multiplication map is $Sp(2g)$-equivariant it follows that
$Sym^2(M_2^\theta(\Gamma_g(2)))\subset M_4^\theta(\Gamma_g(2))$
(so if $f_1,\ldots,f_N$ is a basis of $M_2^\theta(\Gamma_g(2))$
then the $f_if_j$ are linearly independent).

\subsection{Decomposing representations of $Sp(2g)$}\label{casimir}
Given a representation of a finite group on a complex vector 
space, one could determine the value of the character of the representation on each conjugacy class and then use the table of irreducible characters of the group to find the decomposition of
the representation. However, it is very time consuming to compute these character values in our examples. Thus we take another approach, which has the additional advantage of finding explicitly certain subrepresentations. 

There is one conjugacy class of $Sp(2g)$ which has only $2^{2g}-1$
elements, the class of the transvections $t_v$ with $v\in\FF_2^{2g}-\{0\}$, see \ref{transvections}.
If $\rho:Sp(2g)\rightarrow GL(V)$ is a complex representation of $Sp(2g)$, the operator
$$
C\,=\,C_\rho\,:=\,\sum_{v\neq 0}\,\rho(t_v)\qquad(\in GL(V))
$$
obviously satisfies $\rho(g)C\rho(g)^{-1}=C$ for all $g\in Sp(2g)$.
If $V=\oplus V_i^{n_i}$ is the decomposition of $V$ into irreducible 
representations $V_i$, $V_i\not\cong V_j$ if $i\neq j$, then, by Schur's lemma, $C$ must be scalar multiplication by a $\lambda_i\in\CC$
on $V_i$. In particular, the eigenvalues of $C$ are the $\lambda_i$
with multiplicity $n_i(\dim V_i)$ (but it can happen that $\lambda_i=\lambda_j$ for $i\neq j$).

To find $\lambda_i$  we consider the trace of $C$ on $V_i$:
as the $t_v$, $v\neq 0$, are the elements of one conjugacy class, 
$$
Tr(C_{|V_i})\,=\,(2^{2g}-1)Tr(\rho_i(t_v))\,=\,(2^{2g}-1)\chi_i(t_v)
$$
where $v$ is now one specific (but arbitrary) transvection and
$\chi_i$ is the character of the irreducible representation $\rho_i$.
On the other hand,
$$
Tr(C_{|V_i})\,=\,(\dim V_i)\lambda_i,
\qquad\mbox{hence}\quad
\lambda_i\,=\,\frac{(2^{2g}-1)\chi_i(t_v)}{\dim V_i}.
$$
Note that $\ker (C-\lambda I)$ will be the direct sum of the
$V_i^{n_i}$ with $\lambda_i=\lambda$, 
so we do not only get information on the multiplicities 
of the irreducible constituents of $\rho$ but also
on the corresponding subspaces of $V$.

\section{The case $g=1$}\label{g=1}

\subsection{The geometry}
In case $g=1$, the holomorphic map
$$
\HH_1\longrightarrow \PP^1,\qquad\tau\longmapsto (\Theta[0](\tau):\Theta[1](\tau))
$$
is non-constant and thus its image contains an open set of $\PP^1$. This implies
that there are no algebraic relations between the $\Theta[\sigma]$. In fact, the image is $\PP^1-\{6\; \mbox{points}\}$ (these points  correspond to the zeroes of $\eta^{12}$, see below) 
and it is isomorphic to the modular curve $\HH_1/\Gamma_1(2,4)\cong\HH_1/\Gamma_1(4)$.

\subsection{The modular forms of even weight}
The dimension formula from section \ref{dimension formula}
and the fact that $M_{2k}^\theta(\Gamma_1(2))=M_{2k}(\Gamma_1(2))$
shows that $\dim M_{2k}(\Gamma_1(2))=k+1$.
A basis of $M_2(\Gamma_1(2))$ is given by the Heisenberg invariants
$\Theta[0]^4+\Theta[1]^4$ and $(\Theta[0]\Theta[1])^2$. More generally,
a basis of $M_{2k}(\Gamma_1(2))$ is given by the homogeneous polynomials of degree $k$ in these generators: 
$$
(\Theta[0]^4+\Theta[1]^4)^k,\quad
(\Theta[0]\Theta[1])^2(\Theta[0]^4+\Theta[1]^4)^{k-1},\quad\ldots,\quad
(\Theta[0]\Theta[1])^{2k}.
$$ 

\subsection{The $Sp(2)$ representation on $M_{2k}(\Gamma_1(2))$}
The functions $\theta[\Delta]^4$ are in $M_2(\Gamma_1(2))$.
The classical transformation theory of theta functions gives:
{\renewcommand{\arraystretch}{1.2}
$$
\rho_2(S):\left\{\begin{array}{ccr} 
\theta[{}^0_0]^4&\longmapsto& -\theta[{}^0_0]^4\\
\theta[{}^0_1]^4&\longmapsto& -\theta[{}^1_0]^4\\
\theta[{}^1_0]^4&\longmapsto& -\theta[{}^0_1]^4
\end{array}\right.,\qquad
\rho_2(T):
\left\{\begin{array}{ccr} 
\theta[{}^0_0]^4&\longmapsto &\theta[{}^0_1]^4\\
\theta[{}^0_1]^4&\longmapsto &\theta[{}^0_0]^4\\
\theta[{}^1_0]^4&\longmapsto &-\theta[{}^1_0]^4
\end{array}\right.,
$$
}
where $S,T$ are the standard generators of $SL(2,\ZZ)$.
Recall the Jacobi relation $\theta[{}^0_0]^4=\theta[{}^0_1]^4+\theta[{}^1_0]^4$ when computing the matrix of the $\rho_2(g)$'s w.r.t. to the basis of $M_2(\Gamma_1(2))$ given by $\theta[{}^0_0]^4,\theta[{}^0_1]^4$ for example.

In particular,
$M_2(\Gamma_1(2))$ is the unique irreducible two dimensional representation
$\rho_{[21]}$ of $S_3$, hence
$$
M_{2k}(\Gamma_1(2))\,\cong\,Sym^k(\rho_{[21]}).
$$
The group $O^+(2)$ is the group of order two generated by the image of $S\in SL(2,\ZZ)$ in $Sp(2)$.

\subsection{The decomposition of $M_6(\Gamma_1(2))$}
\label{dec g=1}
We briefly discuss the decomposition of $M_6(\Gamma_1(2))$
because it was (implicitly) used in our earlier papers.
As $M_6(\Gamma_1(2))=Sym^3(M_2(\Gamma_1(2)))$, it
is the direct sum of the three irreducible representations of $S_3$:
$$
M_6(\Gamma_1(2))\,\cong\,\rho_{[3]}\oplus\rho_{[2,1]}\oplus\rho_{[1^3]},
$$
where $\rho_{[1^3]} $
is the sign representation of $S_3$ and $\rho_{[3]}$ is the trivial representation.

One verifies that:
{\renewcommand{\arraystretch}{1.5}
$$
\begin{array}{rcl}
\rho_{[3]}\,&=&\,
\langle\,\Theta[0]^{12}-33\Theta[0]^8\Theta[1]^4
-33\Theta[0]^4\Theta[1]^8+\Theta[1]^{12}\rangle,\\
\rho_{[1^3]}\,&=&\,\langle\,\eta^{12}\,\rangle,\\
\rho_{[2,1]}\,&=&\,\{
(a+b)\theta[{}^0_0]^{12}+a \theta[{}^0_1]^{12}+b\theta[{}^1_0]^{12}:\,a,b\in\CC\},
\end{array}
$$
} 
were we used a 
classical formula for the Dedekind $\eta$ function:
$\eta^3=\theta[{}^0_0]\theta[{}^0_1]\theta[{}^1_0]$, so
{\renewcommand{\arraystretch}{1.5}
$$
\begin{array}{rcl}
\eta^{12}&=&\theta[{}^0_0]^4\theta[{}^0_1]^4\theta[{}^1_0]^4\\
&=&(\Theta[0]^2+\Theta[1]^2)^2
(\Theta[0]^2-\Theta[1]^2)^2
(2\Theta[0]\Theta[1])^2.
\end{array}
$$
}
The function $\eta^{12}$ is a modular form of weight 6 on $\Gamma_1(2)$; viewed as homogeneous polynomial of degree $6$, 
it vanishes in the six points $(1:0)$, $(0:1)$, $(1:i^k)$, 
with $k=0,1,2,3$ and $i^2=-1$, in $\PP^1$.

The subspace of $O^+(2)$-anti-invariants is: 
$$
M_6(\Gamma_g(2))^{\epsilon}\,=\,\langle \eta^{12},\,
f_{21}:=2\theta[{}^0_0]^{12}+\theta[{}^0_1]^{12}+\theta[{}^1_0]^{12}
\,\rangle,
$$
the function $f_{21}$ lies in the two-dimensional irreducible subrepresentation.

\section{The case $g=2$}\label{g=2}
\subsection{The geometry}
In case $g=2$, the four $\Theta[\sigma]$'s define a holomorphic map
$$
\HH_2\,\longrightarrow\,\PP^3,\qquad\tau\longmapsto 
(\Theta[00](\tau):\Theta[01](\tau):\Theta[10](\tau):\Theta[11](\tau))
$$
which factors over $\HH_2/\Gamma_2(2,4)$. 
Its image contains an open subset of $\PP^3$, 
hence there are no algebraic relations between the 
$\Theta[\sigma]$'s.
Results of Igusa imply that this map
induces an isomorphism between the Satake compactification of 
$\HH_2/\Gamma_2(2,4)$ and $\PP^3$. 
In particular, the image of the map is $\HH_2/\Gamma_2(2,4)$ and it is $\PP^3$ minus the union of $15$ pairs of lines. 
The pairs of lines are the pairs of eigenspaces of the $15$ elements $(1,x,u)$, with $(x,u)\in \FF_2^4-\{0\}$, 
in the Heisenberg group acting on $\PP^3$. 
For example $(1,(0,0),(1,0))$ maps $\Theta[ab]$ to $(-1)^a\Theta[ab]$ 
and thus has two eigenlines parametrized by $(s:t:0:0)$ and $(0:0:s:t)$ where $(s:t)\in\PP^1$.

The quotient of $\PP^3$ by the Heisenberg group is the Satake compactification of $\HH_2/\Gamma_2(2)$. The ring of Heisenberg invariants can be shown to be generated by $5$ polynomials $p_0,\ldots,p_4$, 
homogeneous of degree 4, in the $\Theta[\sigma]$'s:
$$
\CC[\ldots,\Theta[\sigma],\ldots]^{H_2}\,=\,\CC[p_0,\ldots,p_4]
$$
where the $p_i$ are defined as:
{\renewcommand{\arraystretch}{1.5}
$$
\begin{array}{c}
p_0=\Theta[00]^4+\Theta[01]^4+\Theta[10]^4+\Theta[11]^4,\qquad
p_1=2(\Theta[00]^2\Theta[01]^2+\Theta[10]^2\Theta[11]^2),\\
p_2=2(\Theta[00]^2\Theta[10]^2+\Theta[01]^2\Theta[11]^2),\qquad
p_3=2(\Theta[00]^2\Theta[11]^2+\Theta[01]^2\Theta[10]^2),\\
p_4=4\Theta[00]\Theta[01]\Theta[10]\Theta[11].
\end{array}
$$
}
Thus the quotient map is given by these quartics $\PP^3\rightarrow \PP^4$.
Equivalently, $\HH_2/\Gamma_2(2)$ is the image of the holomorphic map
given by the quartics in the $\Theta[\sigma](\tau)$'s above.
The $15$ pairs of eigenlines map to $15$ lines in $\PP^4$, for example
both $(s:t:0:0)$ and $(0:0:s:t)$ map to
$(s^4+t^4:2s^2t^2:0:0:0)$ which is the general point on the line
parametrized by $(u:v:0:0:0)$ with $(u:v)\in\PP^1$.
The image of $\PP^3$ is defined by a quartic polynomial $f_4$ in 5 variables, the Igusa quartic (cf.\ \cite{CP}), that is $f_4(p_0,\ldots,p_4)$ is identically zero. It follows that, as a graded ring,
$$
\oplus_{k=0}^\infty M_{2k}(\Gamma_2(2))\,\cong\,
\CC[y_0,\ldots,y_4]/(f_4(y_0,\ldots,y_4)).
$$
In particular, $\dim M_2=5$, $\dim M_4=15$, $\dim M_6=35$,
$\dim M_8=70-1=69$.

The singular locus of the image of $\PP^3$, i.e.\ of the threefold in $\PP^4$ defined by $f_4=0$, is the union of the 15 `boundary' lines. 
In \cite{CP} it was shown that the five partial derivatives 
$\partial f_4/\partial y_i(p_0,\ldots,p_4)$, 
span a 5 dimensional space of modular forms of weight $6$
on $\Gamma_2(2)$ which is precisely the space spanned by the $\Xi_6[\delta]$'s, the modular forms introduced by D'Hoker and Phong in 
\cite{DP}.
Hence the $\Xi_6[\delta]$ are zero on the boundary, that is, 
they are cusp forms.

\subsection{The $Sp(4)$-representations on the $M_{2k}(\Gamma_2(2))$}
\label{g=2, k=6}
Using the isomorphism $Sp(4)\cong S_6$, the irreducible representations of $Sp(4)$ are labelled by partitions of $6$. 
One has the identifications
(cf.\ \cite{CP}, one finds $\rho_\theta$ from the representation on $V_\theta:=M_2(\Gamma_2(2))$, $\rho_\theta+\rho_r$ is the representation on the $\theta[\delta]^{12}$, which is $Ind^{Sp}_{O^+}(\epsilon)$
and $1+\sigma_\theta$ is the permutation representation on the $10$ even $\theta[\delta]^8$, hence its trace must be $\geq 0$ on each conjugacy class):
$$
\rho_\theta\,=\,\rho_{[2^3]},\qquad \rho_r\,=\,\rho_{[21^4]},\qquad
\sigma_\theta\,=\,\rho_{[42]}.
$$
Using the representation theory of $S_6$ one finds:
{\renewcommand{\arraystretch}{1.5}
$$
\begin{array}{rccl}
M_2(\Gamma_2(2))\cong&\rho_{[2^3]},&&\\
M_4(\Gamma_2(2))\cong&Sym^2(\rho_{[2^3]})&=&
{\bf 1}+\rho_{[42]}+\rho_{[2^3]},\\
M_6(\Gamma_2(2))\cong&Sym^3(\rho_{[2^3]})&=&
{\bf 1}+2\rho_{[2^3]}+\rho_{[21^4]}+\rho_{[42]}+\rho_{[31^3]},\\
M_8(\Gamma_2(2))\cong&Sym^4(\rho_{[2^3]})\,-\,{\bf 1}&=&
{\bf 1}+3\rho_{[2^3]}+3\rho_{[42]}+\rho_{[31^3]}+\rho_{[321]}.
\end{array}
$$
} 
Explicit functions in the various subrepresentations of $M_6(\Gamma_2(2))$ are given in \cite{CP}.

\section{The case $g=3$}\label{g=3}
\subsection{The geometry} \label{geo g=3}
In case $g=3$, the $8$ $\Theta[\sigma]$'s define a holomorphic map
$$
\HH_3\,\longrightarrow\, \PP^7,\qquad \tau\longmapsto
(\Theta[000](\tau):\ldots:\Theta[111](\tau))
$$
which factors over $\HH_3/\Gamma_3(2,4)$.
The image of this map is a $6$-dimensional quasi-projective variety $Z^0$
whose closure $Z$ is defined by a homogeneous (Heisenberg invariant) polynomial $F_{16}$ of degree $16$, see section \ref{weight 8}. 
In particular, the holomorphic function
$\tau\mapsto F_{16}(\ldots,\Theta[\sigma](\tau),\ldots)$ is identically zero on $\HH_3$. 
The complement of $Z^0$ in $Z$ is the union of 63 pairs
of $\PP^3$'s, which are the eigenspaces of elements in the Heisenberg group.

\subsection{Modular forms}
As there is only one relation of degree $16$ and the quotient ring is normal we get (\cite{R1}, \cite{R2}):
$$
M_{2k}(\Gamma_3(2))\,=\,
M^\theta_{2k}(\Gamma_3(2))\,=\,
(\CC[\ldots,\Theta[\sigma],\ldots]_{4k})^{H_3}
$$ 
with 
{\renewcommand{\arraystretch}{1.5}
$$
(\CC[\ldots,\Theta[\sigma],\ldots]_{4k})^{H_3}\,=
\left\{\begin{array}{ll}
(\CC[\ldots,X_\sigma,\ldots]_{4k})^{H_3}&\quad k\leq 3,\\
(\CC[\ldots,X_\sigma,\ldots]_{4k})^{H_3}/
F_{16}(\CC[\ldots,X_\sigma,\ldots]_{4k-16})^{H_3},&
\quad  k\geq 4.
\end{array}
\right.
$$
} 

\subsection{Weight $2$}\label{weight 2}
From section \ref{weight 2 irr} we know that
$M_2((\Gamma_3(2))$, a fifteen dimensional vector space,
is an irreducible $Sp(6)$-representation, denoted by $\rho_\theta$. 
As $Sp(6)$ has a unique irreducible representation of dimension $15$, denoted by ${\bf 15}_a$ in \cite{Frame}, it follows that
$\rho_\theta\cong {\bf 15}_a$.

\subsection{Weight $4$}\label{weight 4}
From section \ref{sym2 dec} we know that $Sym^2(M_2(\Gamma_3(2)))
\subset M_4(\Gamma_3(2))$
and as an $Sp(6)$-representation we have:
$$
Sym^2(M_2(\Gamma_3(2)))\,:=\,Sym^2({\bf 15}_a)\,=\,
{\bf 1}\,+\,{\bf 35}_b\,+\,{\bf 84}_a.
$$
The invariant subspace is spanned by $\sum_\Delta\theta[\Delta]^8$ and
the subrepresentation ${\bf 1}\,+\,{\bf 35}_b$ is spanned by the $36$
$\theta[\Delta]^8$'s which are permuted by $Sp(6)$.

We recall the following relation, which holds for all $\tau\in\HH_3$:
$$
r_1-r_2=r_3,\qquad {\rm with}\quad
r_1=\prod_{a,b\in\FF_2}\theta[{}^{000}_{0ab}](\tau),\quad
r_2=\prod_{a,b\in\FF_2}\theta[{}^{000}_{1ab}](\tau),\quad
r_3=\prod_{a,b\in\FF_2}\theta[{}^{100}_{0ab}](\tau).
$$
From this we deduce that $2r_1r_2=r_1^2+r_2^2-r_3^2$.
Thus $r_1r_2$, a product of $8$ distinct $\theta[\Delta]$'s,
is a linear combination of three products of four theta squares.
The sum of the four characteristics in each product is zero, hence 
$$
r_1r_2\,=\,\prod_{a,b,c\in\FF_2}\theta[{}^{000}_{abc}]
\,\in M_4(\Gamma_3(2)).
$$
We verified that under the action
of $Sp(6)$ on $r_1r_2$ one obtains $135$ functions which
are a basis of $M_4(\Gamma_3(2))$ and which are permuted (without signs) by $Sp(6)$.

Let $P\subset Sp(6)$ be the stabilizer of $r_1r_2$, it consists
of the matrices with blocks $A,\ldots,D$ with $C=0$.
There are no non-trivial homomorphisms $P\rightarrow GL_1(\CC)=\CC^*$,
because these factor over $SL(3,\FF_2)$ 
(with $P$ as before, one maps the matrix first to $A$) and this is a simple group (of order $168$). 
Thus any $g\in P$ acts as the identity on $r_1r_2$.
By Frobenius reciprocity one can then identify
the representation of $Sp(6)$ on $M_4(\Gamma_3(2))$
with $Ind^{Sp}_P(1_P)$, this representation is (cf.\ \cite{Frame}, p.\ 113)
$$
M_4(\Gamma_3(2))\,\cong\,Ind^{Sp}_P(1_P)\,\cong\, 
{\bf 1}\,+\,{\bf 35}_b\,+\,{\bf 84}_a\,+\,{\bf 15}_a\,\cong\,
Sym^2(M_2(\Gamma_3(2)))\,+\,{\bf 15}_a.
$$

In particular, there is a unique complementary $15$-dimensional
subspace which is $Sp(6)$-invariant; the representation on this
subspace must be ${\bf 15}_a\cong \rho_\theta$. An isomorphism of representations, which uses some geometry of quadrics, can be obtained as follows. 

Let $L\subset \FF_2^6$ be a Lagrangian subspace (so $L\cong \FF_2^3$ and $E(v,w)=0$ for all $v,w\in L$). Then there are eight even quadrics 
$Q$ such that $L\subset Q$ (cf.\ \cite{CDG} Appendix A) In case
$L=L_0=\{(v',v''):v''=0\}$, the 8 even quadrics containing $L_0$ have the
same  characteristics as the eight theta constants in $r_1r_2$.
Each of the $135$ functions in the $Sp(6)$-orbit of $r_1r_2$
can thus be written as 
$$
P_L=\pm \prod_{Q\supset L} \theta[\Delta_Q],
$$
for a unique Lagrangian subspace $L$,
where the product is over the 8 quadrics containing $L$ and the sign is determined by the condition that it is $+1$ if $L=L_0$ and that $P_L=\rho(g)P_{L_0}$ for some $g\in Sp(6)$.

Using this description of the basis of $M_4(\Gamma_3(2))$ it is not hard to write down the (unique) $O^+$-anti-invariant.
Recall that $O^+$ is the stabilizer of the even characteristic $[0]$
and let $Q_0$ be the corresponding even quadric. 
An even (=split) quadric $Q$ in $\FF_2^6$ (cf.\ \ref{quasymp})
contains $30$ Lagrangian subspaces, 
$15$ in each ruling ($L,L'\subset Q$ are in the same ruling if $L\cap L'$ is 1-dimensional). Let $P[0]$ be the sum of the $15$ $P_L$'s from one ruling minus the sum of the  $15$ $P_L$'s from the other ruling
of $Q_0$. Then $P[0]$ transforms with the representation $\epsilon$ of $O^+$. 
Thus the subrepresentation of $M_4(\Gamma_3(2))$ generated by $P[0]$
is contained in $Ind^{Sp}_{O^+}(\epsilon)$ and thus it
must be $\rho_\theta={\bf 15}_a,\,\rho_r={\bf 21}_b$ or their direct sum.
As only ${\bf 15}_a$ is a component of the representation 
on $M_4(\Gamma_3(2))$, we conclude that the subrepresentation generated by $P[0]$ is isomorphic to ${\bf 15}_a$ and thus is complementary to $Sym^2(M_2(\Gamma_3(2)))$.

\subsection{The $Sp(6)$-representation on $M_6(\Gamma_3(2))$}\label{weight 6}
We verified that the $680$ products $P_iP_jP_k$ 
($0\leq i\leq j\leq k\leq 14$) are linearly independent in the $870$-dimensional
vector space $(\CC[\ldots,\Theta[\sigma],\ldots]_{8})^{H_g}$. 
We denote the subspace which they span by
$$
Sym^3(M_2(\Gamma_3(2)))\,:=\,
\langle \,P_iP_jP_k\,:\;P_i,P_j,P_k\in M_2(\Gamma_3(2))\,\rangle.
$$
As $Sp(6,\FF_2)$-representation we have $M_2(\Gamma_3(2))\cong{\bf 15}_a$, and this implies:
$$
Sym^3(M_2(\Gamma_3(2)))\,\cong \,
2\cdot{\bf 15}_a\,+\, {\bf 21}_b\,+\,{\bf 35}_b\,+\,{\bf 84}_a\,+\,
{\bf 105}_c \,+\,{\bf 189}_c + {\bf 216}_a.
$$

To decompose all of $V=M_6(\Gamma_3(2))$ we use operator $C=C_V$ as in
\ref{casimir}. A computer computation showed that the eigenvalues 
$\lambda$ with multiplicities $m_\lambda$ of $C$ on $V$ are:
$$
(\lambda,m_\lambda)\,:\;
(63,1),\quad (27,35),\quad(3,378),\quad(-7,216),\quad(-13,189),
\quad(-21,30),\quad(-33,21).
$$
Obviously $(63,1)$ corresponds to the one-dimensional trivial
representation. 
Using the character table of $Sp(6)$ in \cite{Frame}, p.114--115 (where $t_v$ is in the $16^{\rm th}$ conjugacy class 
labelled $1^{-5}2^6$), one finds unique irreducible representations $\rho$ of dimension $d_\lambda$ such that $C_\rho$ has eigenvalue $\lambda$ in the cases $(\lambda,d_\lambda)=(-7,216)$, $(-13,189)$,
$(-33,21)$. It follows that the irreducible representations 
${\bf 216}_a$, ${\bf 189}_c$ and ${\bf 21}_b$ occur in $M_6(\Gamma_3(2))$, all three with multiplicity one.
The irreducible representations $\rho$ for which $C_\rho$ has eigenvalue $\lambda=27$
are ${\bf 21}_a$ and ${\bf 35}_b$. 
As $C$ has a $35$-dimensional eigenspace with $\lambda=27$, 
the representation
${\bf 35}_b$ occurs with multiplicity one in $M_6(\Gamma_3(2))$.
The irreducible representations $\rho$ for which
$C_\rho$ has eigenvalue $\lambda=3$ are:
${\bf 105}_c$, ${\bf 84}_a$, ${\bf 420}_a$, ${\bf 210}_b$.
This gives two possibilities for the decomposition of the 
$378$-dimensional  eigenspace with $\lambda=3$:
$2\cdot{\bf 105}_c+2\cdot{\bf 84}_a$ or 
${\bf 210}_b+2\cdot{\bf 84}_a$.
As we saw above, ${\bf 105}_c$ is an irreducible component
of $Sym^3({\bf 15}_a)\subset M_6(\Gamma_3(2))$, 
hence we conclude that
$2\cdot{\bf 105}_c+2\cdot{\bf 84}_a$ 
is a summand of $M_6(\Gamma_3(2))$.
The irreducible representations $\rho$ for which
$C_\rho$ has eigenvalue $\lambda=-21$ are:
${\bf 105}_a$ and ${\bf 15}_a$.
As $C$ has a $30$-dimensional eigenspace with $\lambda=-21$ it follows
that $2\cdot {\bf 15}_a$ is the representation on this eigenspace. 

Combining these results with the decomposition of $Sym^3({\bf 15}_a)$ above, we get
$$
M_6(\Gamma_3(2))\,=\,Sym^3({\bf 15}_a)\,+\,{\bf 1}\,+\,{\bf 84}_a\,
+\,{\bf 105}_c.
$$

\subsection{The invariant in $M_6(\Gamma_3(2))$}
The trivial subrepresentation corresponds to an $Sp(6)$-invariant
modular form, that is, to a modular form of weight $6$ for $\Gamma_3$ whose existence follows already from the 
dimension formula for $M_{k}(\Gamma_3)$ found by Tsuyumine \cite{Ts},
c.f.\ \cite{R2}, p.188. 

Our results on the decomposition of $M_{2k}(\Gamma_3(2))$ for $k=1,2$ suggest a way to write down this invariant $G$ which may be of some interest. 
Both $M_{2}(\Gamma_3(2))$ and $M_{4}(\Gamma_3(2))$
have a unique copy of the the $Sp(6)$-representation 
$\rho_\theta={\bf 15}_a$. 
The representation $M_{2}(\Gamma_3(2))\cong \rho_\theta$
contains the $36$ functions $\theta[\Delta]^4$ and in \ref{weight 4} we
identified a function $P[0]$ in 
${\bf 15}_a\subset M_{4}(\Gamma_3(2))$
on which $O^+(6)$ acts via the character $\epsilon$, as it does on $\theta[0]^4$. Since the subspace of $O^+$-anti-invariants in $\rho_\theta$ is one dimensional, we get an isomorphism of $Sp(6)$-representations $\Phi$ by defining
$$ 
\Phi:M_2(\Gamma_3(2))\,\stackrel{\cong}{\longrightarrow}\,
{\bf 15}_a\;(\subset M_4(\Gamma_3(2)))\qquad 
\Phi(\rho_2(g)\theta[0]^4):=\rho_4(g)P[0].
$$
Now let $P[\Delta]:=\Phi(\theta[\Delta]^4)$. Then the function
$G:=\sum_\Delta\theta[\Delta]^4P[\Delta]\in M_6(\Gamma_3(2))$
corresponds to the function 
$\sum_\Delta\theta[\Delta]^8\in Sym^2({\bf 15}_a)$, which is an $Sp(6)$-invariant. Thus $G$ is an $Sp(6)$-invariant and an explicit
verification showed that $G\neq 0$.

\subsection{Asyzygous sextets}\label{as6}
The subrepresentation space ${\bf 105}_c$ in the complement of $Sym^3({\bf 15}_a)$ is quite interesting and,
as we will show, has in fact been studied by D'Hoker and Phong in \cite{DP2}.
Following \cite{DP2} we consider sets 
$S=\{\Delta_1,\ldots,\Delta_6\}$ of six totally asyzygous 
(even) characteristics, that is, 
$\Delta_i+\Delta_j+\Delta_k$ is odd for distinct $i,j,k$.
An example of such a sextet of even characteristics is
$$
S_0\,:=\,\{\;[{}^{110}_{110}],\quad [{}^{110}_{111}],\quad
[{}^{111}_{110}],\quad [{}^{101}_{101}],\quad
[{}^{101}_{111}],\quad [{}^{111}_{101}]\;\}.
$$
Note that the six characteristics in the sextet $S_0$ are of the form 
$[{}^{1ab}_{1cd}]$ where $[{}^{ab}_{cd}]$ runs over 
the six odd theta characteristics in genus $2$.
It is well known that the sum of any three odd characteristics
in genus two is even, so the sum of any three of these six characteristics for $g=3$ is indeed odd. There are $336$ asyzygous
sextets and $Sp(6)$ acts transitively on the $336$ totally asyzygous sextets. The sum of the six characteristics in such a sextet is zero (in $\FF_2$),
hence to a sextet $S$ we can associate a modular form $F_S$
of weight $6$ on $\Gamma_3(2)$:
$$
F_S\,:=\,\prod_{\Delta\in S}\,\theta[\Delta]^2\qquad 
(\in M_6(\Gamma_3(2))).
$$
Using the classical theta formula \ref{classical},
each $F_S$ corresponds to a Heisenberg invariant homogeneous polynomial of degree $12$ in the $\Theta[\sigma]$'s.
A computer computation shows that the $336$ functions $F_S$ span a $105$-dimensional vector space, which we denote by
$$
W_{as}\,:=\,\langle\;F_S\,:\;S\;\mbox{totally asyzygous sextet}\,\rangle
\qquad(\subset M_6(\Gamma_3(2))).
$$
We verified that $Sym^3(M_2(\Gamma_3(2)))\cap W_{as}=\{0\}$ and that
$C$ is multiplication by $3$ on $W_{as}$, hence
$
W_{as}\,\cong\,{\bf 105}_c
$. 
Therefore $W_{as}$ is the 105-dimensional irreducible representation  
in the complement of $Sym^3({\bf 15}_a)$.

One of the results in \cite{DP2} is that there are no functions in
$W_{as}$ which transform as $\theta[0]^4$ under $O^+(6)$.
This is now clear, since such a function would generate a 
subrepresentation of 
$Ind^{Sp}_{O^+}(\epsilon)={\bf 15}_a\oplus{\bf 21}_b$, which contradicts the fact that $W_{as}\cong{\bf 105}_c$.

\subsection{Cusp forms in $W_{as}$}\label{cusp DP}
A modular form $f\in M_{2k}(\Gamma_g(2))$ is called a cusp form if,
for all $M\in \Gamma_g(1)$, and all $\tau_{g-1}\in\HH_{g-1}$ one has,
with $\tau_1\in\HH_1$:
$$
\lim_{\tau_1\rightarrow i\infty}\,
f(M\cdot \tau_{g-1,1})\,=\,0,\qquad
\mbox{with}\quad \tau_{g-1,1}\,=\,
\begin{pmatrix}\tau_{g-1}&0\\0&\tau_1\end{pmatrix}.
$$

We now show that the $F_S$ are cusp forms. The group $Sp(2)$ acts transitively on these functions, so it suffices to show that $F_{S_0}$ 
is a cusp form. 
But $F_{S_0}(M\cdot \tau_{g-1,1})=\pm F_{T}(\tau_{g-1,1})$
where $T$ is sextet $M^{-1}\cdot \Delta_i$, where the $\Delta_i$, $1\leq i\leq 6$ are the characteristics in $S_0$. So we must show that
$\lim_{\tau_1\rightarrow i\infty}\,
F_S(\tau_{g-1,1})=0$ for all $336$ asyzygous sextets $S$.
As the limit $\theta[{}^{abc}_{def}](\tau_{g-1,1})=0$ 
for all $\tau_{g-1}$ iff $c=1$, 
we must verify that each sextet has 
at least one characteristic with $c=1$. 
This can be done by computer, and we did that
(there is a more elegant method, using that the asyzygous sextets correspond to decompositions $\FF_2^{6}=V_1\oplus V_2$ with $V_i$ symplectic subspaces of dimension $V_i$, the condition that $\theta[\Delta]\rightarrow 0$ is just that 
$v=({}^{001}_{000})\in \FF_2^6$ does not lie in $Q_\Delta$ etc.).

\subsection{Miyawaki's $F_{12}$}
The cusp forms $F_S$ are permuted, up to signs, by the 
elements of $Sp(6)$, thus the sum (over the $336$ asyzygous sextets) of their squares 
$$
F_{12}\,=\,\sum_{S}\,F_S^2
$$
is an $Sp(6)$-invariant cusp form of weight $12$, 
if it isn't identically zero on $\HH_3$.
To exclude that possibility we computed the degree $24$ polynomial
in the $X_\sigma$ corresponding to $F_{12}$ and verified that the result
is not a multiple of the polynomial $F_{16}$ (we actually put $X_\sigma=1$ except for $X_{111}$ in this computation). 
Thus $F_{12}$ is a cusp form of weight $12$ on $\Gamma_3(1)$ and was first discovered by I.\ Miyawaki (\cite{Miyawaki}). 
Recently B.\ Heim \cite{Heim} determined the spinor L-function of this cusp form which is of interest for arithmetical applications.

\subsection{The $Sp(6)$-representation on $M_8(\Gamma_3(2))$}
\label{weight 8}
The modular forms of weight $8$ on $\Gamma_3(2)$ can be described as follows:
$$
\CC[\ldots,X_\sigma,\ldots]_{16}^{H_g}\,\cong\,
M_8(\Gamma_3(2)) \oplus \langle F_{16}\rangle,
$$
where $F_{16}$ is a homogeneous polynomial of degree $16$ in the $X_\sigma$ such that 
$F_{16}(\ldots,\Theta[\sigma](\tau),\ldots)=0$ for all $\tau\in\HH_3$,
see \cite{vGvdG}, \cite{CDG}, $\S$ 4.1. This polynomial may be written
as
$$
F_{16}\,=\,8\sum_\Delta\theta[\Delta]^{16}\;-\,
\Bigl(\sum_\Delta \,\theta[\Delta]^{8}\Bigr)^2,\qquad
\theta[{}^\epsilon_{\epsilon'}]^2\,=\,\sum_{\sigma}
(-1)^{\sigma\epsilon'}X_\sigma X_{\sigma+\epsilon}
$$
so one substitutes $X_\sigma$ rather than $\Theta[\sigma]$ in the classical theta formula.
In particular,
$$
\dim M_8(\Gamma_3(2))\,=\,3993-1\,=3992.
$$

We computed the dimension of the image, denoted by $Sym^4(M_2(\Gamma_3(2)))_0$, of  
$Sym^4(M_2(\Gamma_3(2)))$ in this $3992$-dimensional vector space
(that is, of the space spanned by the  products $P_iP_jP_kP_l$, 
for a basis $P_i\in \CC[\ldots,X_\sigma,\ldots]_4^{H_g}$, $1\leq i\leq 15$, the image includes $F_{16}$).
There are $27$ independent relations,
thus
$$
\dim Sym^4(M_2(\Gamma_3(2)))_0=(3060-1)-27=3032.
$$
As $M_2(\Gamma_3(2))={\bf 15}_a$, a computation with $Sp(6)$-representations shows that
$$
\begin{array}{rcl}
Sym^4(M_2(\Gamma_3(2)))\,&\cong&\,
2\cdot{\bf 1}+ 2\cdot{\bf 15}_a+{\bf 27}_a+{\bf 35}_a+4\cdot{\bf 35}_b+
4\cdot {\bf 84}_a+{\bf 105}_c+\\&&+{\bf 168}_a
+{\bf 189}_c+2\cdot{\bf 216}+
3\cdot{\bf 280}_b+{\bf 336}_a+{\bf 420}_a.
\end{array}
$$
The $1+27$-dimensional kernel of the map
$Sym^4(M_2(\Gamma_3(2)))\rightarrow Sym^4(M_2(\Gamma_3(2)))_0$
must be an $Sp(6)$-representation, 
thus it must be ${\bf 1}+{\bf 27}_a$. We have an explicit
description of this kernel and we will show it is useful in the study of modular forms of higher genus in another paper.

It remains to identify the complementary $Sp(6)$-representation, which we will denote by $W$:
$$
M_8(\Gamma_3(2))\,\cong \,Sym^4(M_2(\Gamma_3(2)))_0\,\oplus\,W,
\qquad\dim W=960.
$$

We computed the operator $C$ from section \ref{casimir} on 
$
\CC[\ldots,X_\sigma,\ldots]_{16}^{H_g}.
$
The resulting pairs of the eigenvalues $\lambda$ with multiplicity $m_\lambda$ of $C$ are:
$$
\begin{array}{rrrrrr}
(\lambda,m_\lambda):&(3,1050),&\quad (9,840),&\quad (15,168),&\quad (27,140),&\quad (63,2),\\&
(-3,672),&\quad (-7,648),&\quad (-9,35),&\quad (-13,378),&\quad (-21,60).
\end{array}
$$
Note that this result also implies that ${\bf 27}_a$ cannot be a subrepresentation of $\CC[\ldots,X_\sigma,\ldots]_{16}^{H_g}$ because the eigenvalue $\lambda$ of $C$ on ${\bf 27}_a$ would have been $35$, but this is not an eigenvalue of $C$.

The eigenvalue $63$ corresponds to the subspace of invariants. 
The eigenvalues
$\lambda=9$, $-3$, $-7$, $-13$ occur only on the irreducible representations
${\bf 280}_b$, ${\bf 336}_a$, ${\bf 216}_a$, ${\bf 189}_c$ respectively, hence
$3\cdot {\bf 280}_b+2\cdot{\bf 336}_a +3\cdot {\bf 216}_a+2\cdot {\bf 189}_c$ is a summand of $M_8(\Gamma_3(2))$. 
Thus $W$ has a summand ${\bf 189}_c +{\bf 216}_a+ {\bf 336}_a$.

The eigenvalue $3$ occurs only on ${\bf 84}_a$, ${\bf 105}_c$, 
${\bf 210}_b$ and  ${\bf 420}_a$. As $4\cdot {\bf 84}_a+{\bf 105}_c
+{\bf 420}_a$ is a summand of $Sym^4(M_2(\Gamma_3(2)))_0$, there
remains a subrepresentation of dimension $1050-861=189$ in $W$ with
this eigenvalue. That implies that ${\bf 84}_a+{\bf 105}_c$
is a summand of $W$. 

The eigenvalue $15$ occurs only on ${\bf 105}_b$, ${\bf 168}_a$
and  ${\bf 210}_a$, 
as the eigenspace of $C$ for this eigenvalue has dimension $168$
we conclude that ${\bf 168}_a$ is a summand of $M_8(\Gamma_3(2))$
(which lies in $Sym^4({\bf 15}_a)$).

The eigenvalue $27$ occurs only on ${\bf 21}_a$ and ${\bf 35}_b$.
As $4\cdot {\bf 35}_b$ is a summand of $Sym^4(M_2(\Gamma_3(2)))_0$
and the dimension of this eigenspace of $C$ is $140$, none of these two representations occurs in $W$.

The eigenvalue $-9$ occurs on four irreducible representations, 
the one with lowest dimension is ${\bf 35}_a$. 
As the eigenspace for $\lambda=-9$ has dimension $35$, 
we conclude that ${\bf 35}_a$ is a summand of $M_8(\Gamma_3(2))$
(which lies in $Sym^4({\bf 15}_a)$).

The eigenvalue $-21$ occurs only on ${\bf 15}_a$ and ${\bf 105}_a$,
as the eigenspace of $C$ for this eigenvalue has dimension $60$
we conclude that $4\cdot  {\bf 15}_a$ is a summand of $M_8(\Gamma_3(2))$
and that $W$ has a summand $2\cdot{\bf 15}_a$.

From this we find the following decomposition of $W$:
$$
W\,=\,2\cdot{\bf 15}_a+{\bf 84}_a+{\bf 105}_c+{\bf 189}_c+
{\bf 216}_a+{\bf 336}_a.
$$

\section{The case $g=4$}\label{g=4}
\subsection{The $Sp(8)$-representation on $M^\theta_{4}(\Gamma_4(2))$}\label{g=4, k=2}
In case $g=4$ we are no longer sure if $M^\theta_{2k}(\Gamma_4(2))$,
the space of modular forms of weight of $2k$ which are (Heisenberg-invariant) polynomials in the $\Theta[\sigma]$'s, is equal
to all of $M_{2k}(\Gamma_4(2))$ (cf.\ \cite{OSM}). 

The $Sp(8)$-representation on $M^\theta_2(\Gamma_4(2))$,
which we denoted by $\rho_\theta$ in section \ref{weight 2 irr},
is the unique $51$-dimensional irreducible representation of $Sp(8)$
(a table of the 
$81$ irreducible representations of $Sp(8)$ can be easily generated with the computer algebra program `Magma' \cite{magma}).

The complement of $Sym^2(M^\theta_2(\Gamma_4(2))$ in
$M^\theta_4(\Gamma_4(2))$ (of $Sp(8)$-representations)
now has codimension $918$:
$$
\dim M^\theta_4(\Gamma_4(2))\,-\,\dim Sym^2(M^\theta_2(\Gamma_4(2)))
\,=\,
2244-\binom{51+1}{2}\,=\,2244-1326\,=\,918.
$$

We computed the operator $C$ from section \ref{casimir} on 
$M_4^\theta(\Gamma_4(2))$. The resulting pairs of its eigenvalues $\lambda$ with multiplicity $m_\lambda$ are:
$$
(\lambda,m_\lambda):\qquad (-25,918),\quad (39,1190),\quad (119,135), \quad (255,1).
$$
The last three eigenspaces of $C$ correspond to the irreducible representations $\sigma_c$, $\sigma_\theta$ and ${\bf 1}$ respectively
and their direct sum is $Sym^2(\rho_\theta)$, cf.\ \ref{sym2 dec}. 
The character table shows that there are only $10$ irreducible representations of $Sp(8)$ with
dimension less then $918$ and there is a unique irreducible representation
with dimension $918$. However, of these eleven irreducible
representations, the map $C$ has eigenvalue $\lambda=-25$ only on
the one of dimension $918$. 
Thus we conclude that 
$M_4^\theta(\Gamma_4(2))$ is the sum of just four irreducible representations (like $M_4^\theta(\Gamma_3(2))$, 
cf.\ section \ref{weight 4}).

\subsection{Frame's observations on $Sp(2g)$-representations}
Frame (\cite{Flast}) observed that the irreducible representations of $Sp(2g)$ seem to come
in series, parametrized by $g$ and an irreducible representation
of $O^+(2l)$ or $O^-(2l)$ for some $l< g$, he called the integer 
$l$ the level of the representation. 

The dimension of an irreducible representation of level
$l$, as function of $g$, 
should be a polynomial in $2^g$ of degree $2l$.
For example, the representations $\rho_\theta$ (of dimension 
$(2^g+1)(2^{g-1}+1)/3=(2^g+1)(2^{g}+2)/(3!)$) is of level one
and should correspond to a representation of 
$O^-(2)\cong S_3$.

Note that we have:
$$
\dim M^\theta_4(\Gamma_g(2))\,-\,\dim Sym^2(M^\theta_2(\Gamma_g(2)))\,=\,
\frac{1}{8!}(2^g+1)(2^g-1)(2^g+2)(2^g-2)(2^g-2^2)(2^g+2^5)
$$
which suggests, in this frame work, that $M^\theta_4(\Gamma_g(2))$
is the direct sum of 
$Sym^2(M^\theta_2(\Gamma_g(2)))$ (itself a direct sum 
of three irreducible representations with multiplicity one)
and an irreducible representation of level $3$ corresponding to
a representation of $O^+(6)\cong S_8$, generalizing the results in \ref{weight 4}, \ref{g=4, k=2}.

\section{Applications to superstring measures}\label{DP}

\subsection{} In this section we consider applications of our results
to modular forms $\Xi_6[0^{(g)}]$ and $\Xi_8[0^{(g)}]$ on 
$\Gamma_g(2)$, of weight $6$ and $8$ respectively, which appear to play an important role in the theory of superstrings 
(cf. \cite{DP}, \cite{DP1}, \cite{DP2}, \cite{CDG}, \cite{Mo}). 
One has $\Xi_6[0^{(1)}]=\eta^{12}$ and
$\Xi_8[0^{(1)}]=\theta[{}^0_0]^4\eta^{12}$.
In case $g=2$ D'Hoker and Phong determined, using superstring theory, a modular form $\Xi_6[0^{(2)}]\in M_6(\Gamma_2(2))^\epsilon$. 

As it is (still) difficult to extend their methods to higher genera, they stated certain constraints which should be satisfied by the $\Xi_6[0^{(g)}]$.
However, in genus three they could not find a modular form satisfying their constraints and we show in \ref{non xi6g3} that indeed such a 
modular form does not exist. In \cite{CDG} the constraints were modified
and now one has to find $O^+$-invariant modular forms $\Xi_8[0^{(g)}]$
of weight $8$ on $\Gamma_g(2)$ which moreover restrict to 
$\Xi_8[0^{(k)}]\Xi_8[0^{(g-k)}]$ on the subvarieties $\HH_k\times\HH_{g-k}\subset\HH_g$. Actually, superstring theory only requires the $\Xi_8$'s to be defined on the Jacobi locus $J_g$ of period matrices of Riemann surfaces. 

In \cite{CDG}, \cite{CDG2} modular forms
$\Xi_8[0^{(3)}]$ and $\Xi_8[0^{(4)}]$ satisfying the modified constraints were found.  Grushevsky \cite{Grr}
gives a more general approach for all genera but it is not clear that his proposal leads to single valued functions for $g\geq 6$, for $g=5$ see \cite{SM}. 
Here we give uniqueness results for these functions.  
For example in genus two the function $\Xi_6[0^{(2)}]$, found with hard work by D'Hoker and Phong, is easily recovered as the only modular form satisfying their constraints (but the validity in superstring theory of the (modified) constraints is not clear yet).

\subsection{The uniqueness of $\Xi_8[0^{(2)}]$}\label{unique xi8g2}
In \cite{CDG}, section $\S$3.4, we announced that a certain modular form
$\Xi_8[0^{(2)}]\in M_8(\Gamma_2(2))^{O^+}$, that is, a modular form on $\Gamma_g(1,2)$, is the unique modular form satisfying certain constraints. To verify this, we recall
from \ref{dim O+ inv} that $\dim M_8(\Gamma_2(2))^{O^+}=n_1+n_\theta$, with $n_1$, $n_\theta$ the multiplicity of 
${\bf 1}$ and $\sigma_\theta=\rho_{[42]}$ in $M_8(\Gamma_2(2))$ respectively, see section \ref{g=2, k=6}.
From the decomposition  given there we find that 
$$
\dim M_8(\Gamma_2(2))^{O^+}\,\,=\,1+3\,=\,4.
$$ 

The subspace $M_8(\Gamma_2(2))^{O^+}$ contains the $Sp(4)$-invariant 
$\sum_\delta\theta[\delta]^{16}\in M_8(\Gamma_2(2))$ as well as the three dimensional subspace spanned by
$$
f_1=\theta[{}^{00}_{00}]^{16},\quad
f_2=\theta[{}^{00}_{00}]^4\sum_\delta\theta[\delta]^{12},\quad
f_3=\theta[{}^{00}_{00}]^8\sum_\delta\theta[\delta]^{8}.
$$
We checked, using the classical theta formulas, 
that these four functions are linearly independent and thus are a basis 
of $M_8(\Gamma_2(2))^{O^+}$. 

The function $\Xi_8[0^{(2)}]$ should restrict to
$\Xi_8[0^{(1)}](\tau_1)\Xi_8[0^{(1)}](\tau_2)$ with $\Xi_8[0^{(1)}](\tau_1)=(\theta[{}^0_0]^4\eta^{12})(\tau_1)$
on $\HH_1\times\HH_1\subset\HH_2$. This function is a multiple of $\theta[{}^0_0]^4(\tau_1)$. The restrictions of the $f_i$ are also multiples of $\theta[{}^0_0]^4(\tau_1)$, but the restriction of $\sum\theta[\delta]^{16}$ is not. Hence $\Xi_8[0^{(2)}]$
should be linear combination of the three $f_i$, 
and in \cite{CDG} we showed that there is a unique such linear combination satisfying all
the constraints. Hence $\Xi_8[0^{(2)}]$ is unique. Moreover, the 
formula given in \cite{CDG} shows that 
$\Xi_8[0^{(2)}]=\theta[0^{(2)}]^4\Xi_6[0^{(2)}]$ where $\Xi_6[0^{(2)}]$
is the modular form found by D'Hoker and Phong.

\subsection{The non-existence of $\Xi_6[0^{(3)}]$}\label{non xi6g3}
In \cite{DP2} the existence of a modular form $\Xi_6[0^{(3)}]\in
M_6(\Gamma_g(2))^{\epsilon}$ with certain properties is investigated. 
The only $Sp(6)$-representations with $O^+$-anti-invariants
are $\rho_\theta={\bf 15}_a$ and $\rho_r={\bf 21}_b$ 
and these have a unique such 
anti-invariant (cf.\ \ref{dim O+ inv}). 
Thus from the decomposition of $M_6(\Gamma_3(2))$
given in \ref{weight 6} it follows that $\dim M_6(\Gamma_3(2))^{\epsilon}=3$. 
We verified that the following functions
are a basis:
$$
M_6(\Gamma_g(2))^{\epsilon}\,=\,\langle
\;\;\theta[{}^{000}_{000}]^{12},\quad
\sum_\Delta \,\theta[\Delta]^{12},\quad
\theta[{}^{000}_{000}]^4\sum_\Delta \,\theta[\Delta]^{8}\;\rangle.
$$ 
The function $\Xi_6[0^{(3)}]$ should restrict to 
$\Xi_6[0^{(1)}](\tau_1)\Xi_6[0^{(2)}](\tau_2)$ for 
$(\tau_1,\tau_2)\in\HH_1\times\HH_2\subset\HH_3$, with
$\Xi_6[0^{(1)}](\tau_1)=\eta^{12}(\tau_1)$.
The restriction map on the theta constants is given by
$$
\theta[{}^{abc}_{def}](\tau)\,\longmapsto 
\theta[{}^{a}_{d}](\tau_1)\theta[{}^{bc}_{ef}](\tau_2),
$$
in particular $\theta[{}^{abc}_{def}]\mapsto0$ if $ad=1$.
Thus $6$ of the $36$ even theta constants map to zero, the other
$30=3\cdot 10$ are uniquely decomposed in the product of two even theta constants for $g=1$ and $g=2$ respectively.
Hence, using the results from \ref{dec g=1} and \cite{CDG}, we get
{\renewcommand{\arraystretch}{1.5}
$$
\begin{array}{rcl}
\theta[{}^{000}_{000}]^{12}(\tau_{1,2})&=&
\theta[{}^0_0]^{12}(\tau_{1})
\theta[{}^{00}_{00}]^{12}(\tau_{2})\\
&=&(\mbox{$\frac{1}{3}$}f_{21}+\eta^{12})\theta[{}^{00}_{00}]^{12},
\\&&\\
(\sum_\Delta\theta[\Delta]^{12})(\tau_{1,2})&=&
(\theta[{}^0_0]^{12}+\theta[{}^0_1]^{12}+\theta[{}^1_0]^{12})(\tau_1)
(\sum_\delta \theta[\delta]^{12})(\tau_2)\\
&=&(\mbox{$\frac{2}{3}$}f_{21}-\eta^{12})\sum_\delta \theta[\delta]^{12},
\\&&\\
(\theta[{}^{000}_{000}]^4\Psi_4)(\tau_{1,2})&=&
\theta[{}^{000}_{000}]^4(\tau_{1,2})
(\sum_\Delta\theta[\Delta]^8)(\tau_{1,2})\\
&=&\theta[{}^0_0]^4(\tau_{1})
(\theta[{}^0_0]^8+\theta[{}^0_1]^8+\theta[{}^1_0]^8)(\tau_1)
\theta[{}^{00}_{00}]^4(\tau_{2})
(\sum_\delta \theta[\delta]^{8})(\tau_2)\\
&=&\mbox{$\frac{2}{3}$}f_{21}\theta[{}^{00}_{00}]^4(\sum_\delta \theta[\delta]^{8}).\\
&&
\end{array}
$$
} 
Thus $a\theta[{}^{000}_{000}]^{12}+
b\sum_\Delta\theta[\Delta]^{12}+c\theta[{}^{000}_{000}]^4(\sum_\Delta\theta[\Delta]^{8})$ restricts to a function of the form 
$\eta^{12}(\tau_1)g(\tau_2)$ iff
$$
a\theta[{}^{00}_{00}]^{12}\,+\,2b\sum_\delta \theta[\delta]^{12}\,+\,
2c\theta[{}^{00}_{00}]^4\sum_\delta \theta[\delta]^{8}=0
$$
on $\HH_2$. However, it follows from the results in Table 6 of \cite{CP}
that the ten even $\theta[\delta]^{12}$'s 
span a ten dimensional space and that $\theta[{}^{00}_{00}]^4\sum_\delta \theta[\delta]^{8}$ does not lie in that $10$ dimensional space. 
Hence we must have $a=b=c=0$, which proves that there is no function
$\Xi_6[0^{(3)}]$ satisfying the three constraints imposed in \cite{DP2}.

\subsection{The uniqueness of $\Xi_8[0^{(3)}]$}\label{unique xi8g3}
In \cite{CDG} we found a modular form $\Xi_8[{}^{000}_{000}]$ which satisfied constraints similar to those of \cite{DP2} and 
which are expected to be related to the superstring measure.

The constraints imply that $\Xi_8[0^{(3)}]$ is a modular 
form  on $\Gamma_3(2)$ of weight $8$ and should be $O^+$-invariant (equivalently, it is a modular form on $\Gamma_3(1,2)$), so it must lie in $M_8(\Gamma_3(2))^{O^+}$.
The only $Sp(6)$-representations that have $O^+(6)$-invariants
are ${\bf 1}$ and $\sigma_\theta={\bf 35}_b$.
Hence (cf.\ \ref{dim O+ inv}), 
the decomposition of $M_8(\Gamma_3(2))$ obtained in section
\ref{weight 8} implies that
$$
\dim M_8(\Gamma_3(2))^{O^+}\,\,=\,1+4\,=\,5.
$$ 

The subspace $M_8(\Gamma_3(2))^{O^+}$ contains the $Sp(6)$-invariant 
$\sum_\Delta\theta[\Delta]^{16}\in M_8(\Gamma_3(2))$ as well as the three dimensional subspace
$$
\theta[{}^{000}_{000}]^4M_6(\Gamma_3(2))^{\epsilon}
\,:=\,\{\theta[{}^{000}_{000}]^4f:
\;f\in M_6(\Gamma_3(2))^{\epsilon}\,\},
$$
since both $\theta[{}^{000}_{000}]^4$ and such $f$ are $O^+$-anti-invariant.

For a Lagrangian subspace $L\subset \FF_2^6$ 
we defined modular forms $\pm P_L\in M_4(\Gamma_3(2))$ 
(cf.\ \ref{weight 4}).
For an even characteristic $\Delta$, the sum of the $30$ $P_L^2$'s
such that $L\subset Q_\Delta$, where 
$Q_\Delta$ is the quadric in $\FF_2^6$ defined by $\Delta$ 
is a modular form of weight $8$ on $\Gamma_g(2)$:
$$
G[\Delta]\,:=\,\sum_{L\subset Q_\Delta}\,P_L^2\,=\,
\,\sum_{L\subset Q_\Delta}\,
\prod_{Q'\supset L}\theta[\Delta_{Q'}]^2\qquad(\in M_8(\Gamma_3(2))).
$$
Note that $\theta[\Delta]^2$ is one of the factors in each of the 30 products. In \cite{CDG} we verified that $G[0]\in M_8(\Gamma_3(2))^{O^+}$.

It is not hard to verify (using the classical theta formulas for example) that we now have functions which span the $O^+$-invariants:
$$
M_8(\Gamma_8(2))^{O^+}\,=\,
\theta[{}^{000}_{000}]^4M_6(\Gamma_8(2))^{\epsilon}\,\oplus\,
\langle\; \sum_\Delta \theta[\Delta]^{16},\;G[0]\;\rangle.
$$

The modular form $\Xi_8[0^{(0)}]$ 
should moreover restrict to
$(\theta[{}^0_0]^4\eta^{12})(\tau_1)
(\theta[{}^{00}_{00}]^4\Xi_6[{}^{00}_{00}])(\tau_2)$
on $\HH_1\times\HH_2\subset\HH_3$, 
note this restriction is a multiple of $\theta[{}^0_0](\tau_1)$.

The restriction of $\sum \theta[\Delta]^{16}$ to the diagonal is not a multiple of $\theta[{}^0_0](\tau_1)$, whereas the other four functions 
do have a restriction which is a multiple of $\theta[{}^0_0](\tau_1)$.
Thus $\Xi_8[0^{(3)}]$ should be linear combination
$$
\Xi_{a,b,c,d}:=
\theta[{}^{000}_{000}]^4(a\theta[{}^{000}_{000}]^{12}+
b\sum_\Delta\theta[\Delta]^{12}+c\theta[{}^{000}_{000}]^4\Psi_4)
\,+\,dG[{}^{000}_{000}].
$$
An explicit computation (cf.\ \cite{CDG}) shows
that only
$
\Xi_8[0^{(3)}]:=\Xi_{4,4,-3,-12}/12
$
restricts to $(\theta[{}^0_0]^4\eta^{12})(\tau_1)
(\theta[{}^{00}_{00}]^4\Xi_6[{}^{00}_{00}])(\tau_2)$. This verifies the uniqueness of $\Xi_8[0^{(3)}]$.

\subsection{The uniqueness of $\Xi_8[0^{(4)}]$}\label{unique xi8g4}
In \cite{CDG2}, \cite{Grr} a modular form $\Xi_8[0^{(4)}]\in M_8(\Gamma_4(2))^{O^+}$
was defined which restricts to $\Xi_8[0^{(a)}](\tau_a)\Xi_8[0^{(b)}](\tau_b)$ for 
$(\tau_a,\tau_b)\in\HH_a\times\HH_b\subset\HH_4$ 
for $(a,b)=(1,3),(2,2)$. We will assume from now on that $\Xi_8[0^{(4)}]$ actually is in $M_8^\theta(\Gamma_4(2))^{O^+}$ (this is currently under investigation).

Using the same methods as in \ref{unique xi8g2}, \ref{unique xi8g3},
one can show that these properties characterize $\Xi_8[0^{(4)}]$, 
up to adding a term $\lambda J$, 
where $J=0$ defines the Jacobi locus $J_4$ 
(the locus of period matrices of Riemann surfaces in $\HH_4$) and $\lambda\in\CC$. In fact, $J$ is a modular form of weight $8$ on $\Gamma_4$, so it can be added to $\Xi_8[0^{(4)}]$
without changing its $\Gamma_g(1,2)$-invariance, moreover $\HH_1\times\HH_3$ and $\HH_2\times\HH_2$ are contained in the closure of the Jacobi locus, so $J$ is zero on these loci.

The key point is the determination of the dimension of $M_8^\theta(\Gamma_4(2))^{O^+}$, which could be done by computer.
Recently M.\ Oura \cite{Oura} determined this dimension using the methods from \cite{R1}, \cite{R2}:
$\dim M^\theta_8(\Gamma_4(2))^{O^+}=7$. As we already know $7$ independent functions
in $M_8^\theta(\Gamma_4(2))^{O^+}$, it follows
(using the notation form \cite{CDG2}):
$$
M_8^{\theta}(\Gamma_4(2))^{O^+}\,=\,
\langle 
\sum_\Delta\theta[\Delta]^{16},\; 
(\sum_\Delta \theta[\Delta]^8)^2,\;
\theta[0^{(4)}]\tilde{F}_1,\;\theta[0^{(4)}]\tilde{F}_2,\;
\theta[0^{(4)}]\tilde{F}_3,\;
 G_1[0^{(4)}],\;G_2[0^{(4)}]\,\rangle.
$$
The modular form $J$ vanishing on the Jacobi locus is 
$$
J\,=\,16\sum\theta[\Delta]^{16}-(\sum \theta[\Delta]^8)^2
$$ 
(cf.\ \cite{IgusaC}).
Now the proof of the unicity of $\Xi_8[0^{(4)}]$ in $M_8^\theta(\Gamma_g(2))$ can be obtained
with arguments similar to those in sections  \ref{unique xi8g2}, \ref{unique xi8g3}.

\subsection{The cosmological constant in $g=3$, $g=4$}
\label{cosmog3}
In \cite{CDG} $\S$ 4.5 we showed that the cosmological constant, which is the sum over the $36$ even characteristics
of $\Xi_8[\Delta^{(3)}]$, is zero. We used that there is a unique, up to scalar multiple, $Sp(6)$-invariant in $M_8(\Gamma_3(2))$.
This follows directly from the fact that in the decomposition of 
$M_8(\Gamma_3(2))$ the representation ${\bf 1}$ has multiplicity one.

In case $g=4$, Grushevsky \cite{Grr} already proved that the cosmological constant vanishes using results on the slope of divisors on the moduli space of genus four curves. Another proof was given by 
R.\ Salvati Manni in \cite{SM}.
We sketch an alternative approach to this result, 
similar to the $g=3$ case in \cite{CDG}. 
A computer computation showed that the space of $Sp(8)$-invariants in $\CC[\ldots,X_\sigma,\ldots]^{H_g}$ is two dimensional (one starts with finding invariants for the transvections which act `diagonally' on the $X_\sigma$'s to reduce the computation to a manageable size).
These invariants correspond to the modular forms $\Psi_8,\Psi_4^2$
with $\Psi_{4k}(\tau):=\sum\theta[\Delta^{(4)}]^{8k}$, note that $J= 16\Psi_8-\Psi_4^2$ is zero on the Jacobi locus.

As $\sum_\Delta \Xi_8[\Delta^{(4)}]$ is an $Sp(8)$-invariant of weight 8,
there are constants $\lambda,\mu$ such that
$$
\sum_\Delta \Xi_8[\Delta^{(4)}]\,=\,\lambda \Psi_8\,+\,\mu \Psi_4^2
$$
and it suffices to show that $\lambda+16\mu=0$. For this one can specialize $\tau$ to a diagonal matrix $\tau={\rm diag}(\tau_1,\ldots,\tau_4)$ and then let $ \tau_j\mapsto\infty$ for $j=1,\ldots,4$, similar to the computation in \cite{CDG}, $\S$4.5.

\subsection*{Acknowledgments}
We would like to thank Sergio L.\ Cacciatori for several interesting and stimulating discussions, 
R.\ Salvati Manni for useful discussions and references and to
L.\ Paredi for supporting us in computer computations.

\

\begin{appendix}


\section{A: The action of $Sp(2g,\FF_2)$ on the Heisenberg invariants}
\label{A}

\subsection{Characteristics and quadrics}
\label{quasymp}
The group $Sp(2g)$ fixes a non-degenerate alternating form $E$ (i.e.\ $E(v,v)=0$ for
all $v$) on the $\FF_2$-vector space $V=\FF_2^{2g}$.
We choose a symplectic basis of $V$ so that 
$$
E:V\times V\, \longrightarrow\, \FF_2,\qquad
E(v,w):=v_1w_{g+1}+\ldots+v_gw_{2g}+
v_{g+1}w_1+\ldots+v_{2g}w_g,
$$
or, more compactly, $E((v',v''),(w',w''))=\,{}^t\!v'w''+\,{}^t\!v''w'$ and we occasionally write
$v=({}^{v'}_{v''})$, where $v',v''$ are then considered as row vectors.

We consider the quadratic forms $q$ on $V$ whose associated bilinear form is $E$
(see also \cite{CDG} Appendix A, \cite{Igusa}, V.6), that is the maps
$$
q:\,V\,\longrightarrow\,\FF_2,\qquad q(v+w)=q(v)+q(w)+E(v,w).
$$
One verifies easily that for all $a_i,b_i\in\FF_2$
the function
$$
q(v)=v_1v_{g+1}+v_2v_{g+2}+\ldots+v_gv_{2g}+
a_1v_1+\ldots+a_gv_g+b_1v_{g+1}+\ldots+b_gv_{2g}
$$
satisfies $q(v+w)=q(v)+q(w)+E(v,w)$ and that any quadratic form associated to $E$ is of this form.
Note that $q(v)=\,{}^t\!v'v''+av'+bv''$, 
with row vectors $a=(a_1,\ldots,a_g)$, $b=(b_1,\ldots,b_g)$.
The (theta) characteristic $\Delta_q$ associated to $q$ is defined as
$$
\Delta_q:=\left[{}^{a_1\,a_2\,\ldots\,a_g}_
{b_1\,b_2\,\ldots\,b_g}\right]
\,=\,[{}^a_{b}],
\qquad \mbox{let}\quad 
e(\Delta_q):=(-1)^{\sum_{i=1}^g a_ib_i}\quad(\in\{1,-1\}).
$$
We say that $\Delta_q$ is even if $e(\Delta_q)=+1$ and odd else.
One can verify that
$q(v)$ has $2^{g-1}(2^g+1)$ zeroes in $V$ if 
$\Delta_q$ is even  and has $2^{g-1}(2^g-1)$ zeroes if $\Delta_q$ is
odd.
Moreover, there are $2^{g-1}(2^g+1)$ even characteristics and 
$2^{g-1}(2^g-1)$ odd characteristics, thus there are $3,10,36,136$ even characteristics and $1,6,28,120$ odd characteristics for $g=1,2,3,4$ respectively.

The group $Sp(2g)$ acts naturally on the characteristics by
$$
(g\cdot q)(v):=q(g^{-1}v)\qquad(g\in Sp(2g),\;v\in V).
$$
This action is transitive (actually doubly transitive 
\cite{Igusa} V.6 Corollary) on both the set of even characteristic and on the set of odd characteristics. 

\subsection{Theta constants with characteristics}\label{thetas}
For an even characteristic $\Delta=[{}^a_{b}]$ one defines a function, 
a theta constant, on the Siegel space $\HH_g$ by
$$
\theta[{}^a_b](\tau)\,:=\,\sum_{m\in\ZZ^g}
\,e^{\pi i((m+a/2)\tau\,{}^t\!(m+a/2)+(m+a/2)\,{}^t\!b)}
$$
so $m$ is a row vector and $\sum a_ib_i\equiv 0$ mod $2$. Then one has (\cite{Igusa}, V.1, Corollary), with $\kappa(M)e^{2\pi i\phi_{\Delta}(M)}$ an eight-root of unity, that
$$
\theta[M\cdot\Delta](M\cdot\tau)\,=\,
\kappa(M)e^{2\pi i\phi_{\Delta}(M)}\det(C\tau+D)^{1/2}\theta[\Delta](\tau),
$$
for all $M\in Sp(2g,\ZZ)$, here the action of $M$ on the characteristic $\Delta$ is given by
$$
\begin{pmatrix}A&B\\C&D\end{pmatrix}\cdot [{}^a_b]\,:=\,
[{}^c_d],\qquad
\left(\begin{array}{c}{}^tc\\{}^td\end{array}\right)\,=\,
\left(\begin{array}{cc} D&-C\\-B&A\end{array}\right)
\left(\begin{array}{c}{}^ta\\{}^tb\end{array}\right)
\,+\,
\left(\begin{array}{c}{}^t(C\,{}^t\!\!D)_0\\{}^t(A\,{}^t\!\!B)_0\end{array}\right)
\quad \mbox{mod}\;2
$$
where $N_0=(N_{11},\ldots,N_{gg})$ is the row vector of diagonal 
entries of the matrix $N$.

To get representations of (quotients of) $\Gamma_g$, we defined the action of
$\Gamma_g$ on modular forms of weight $k$ by 
$(M\cdot f)(\tau):=(\det(C\tau+D)^k f(M^{-1}\cdot\tau)$ (cf. \ref{level 2}).
Note that
$$
\theta[\Delta](M^{-1}\cdot\tau)\,=\,
\theta[M^{-1}M\cdot\Delta](M^{-1}\cdot\tau)\,=\,
c_{M^{-1},\Delta,\tau}\theta[M\cdot\Delta](\tau)
$$
where $c_{M^{-1},\Delta,\tau}$ collects the non-relevant part.
Thus the action of $M$ basically maps $\theta[\Delta]$ to
$\theta[M\cdot\Delta]$. 

\subsection{The action of $Sp(2g)$ on the characteristics}
\label{sp char}
We show explicitly that 
$$
(M\cdot q_\Delta)(v)\,=\,q_{M\cdot \Delta}(v),
$$
which verifies that the natural action of $\Gamma_g$ on the theta characteristics corresponds to its action on the quadratic forms on $V$.

As $
(M\cdot q_\Delta)(v)\,=\,q_\Delta(M^{-1}v)$, we must verify 
$q_\Delta(M^{-1}v)=q_{M\cdot \Delta}(v)$.
For $M\in\Gamma_g=Sp(2g,\ZZ)$ we have
$M E\,{}^t\!M =E$ where $E$ has blocks $A=D=0,B=-C=I$:
$$
M E\,{}^t\!M =E\qquad\mbox{iff}\quad
\begin{pmatrix}-B\,{}^t\!A+A\,{}^t\!B&-B\,{}^t\!C+A\,{}^t\!D\\
-D\,{}^t\!A+C\,{}^t\!B&-D\,{}^t\!C+C\,{}^t\!D\end{pmatrix}=
\begin{pmatrix}0&I\\-I&0\end{pmatrix}.
$$
As $E^{-1}=-E$ we find
that
$$
M^{-1}=-E\,{}^t\!M E,\qquad{\rm so}\quad
 M^{-1}v=\begin{pmatrix}\,{}^t\!D&-\,{}^t\!B\\-\,{}^t\!C&\,{}^t\!A\end{pmatrix}
\begin{pmatrix}v'\\v''\end{pmatrix}\,=\,
\begin{pmatrix}\,{}^t\!Dv'-\,{}^t\!Bv''\\-\,{}^t\!Cv'+\,{}^t\!Av''\end{pmatrix}.
$$
Therefore, with $q_\Delta(v)\,=\,q_{[{}^a_b]}(({}^{v'}_{v''}))\,=\,
\,{}^t\!v'v''+av'+bv''$, we get:
$$
q_\Delta(M^{-1}v)\,=\,\,{}^t\!(\,{}^t\!Dv'-\,{}^t\!Bv'')(-\,{}^t\!Cv'+\,{}^t\!Av'')+
a(\,{}^t\!Dv'-\,{}^t\!Bv'')+b(-\,{}^t\!Cv'+\,{}^t\!Av'')
$$
The non-linear part is
$$
\,{}^t\!v'(-D\,{}^t\!C)v'+\,{}^t\!v'(D\,{}^t\!A)v''+
\,{}^t\!v''(B\,{}^t\!C)v'-\,{}^t\!v''(B\,{}^t\!A)v''.
$$
As $M$ is symplectic, $D\,{}^t\!C$ is symmetric, hence only the terms $(D\,{}^t\!C)_{ii}(v'_i)^2$ remain mod 2, but $(v'_i)^2\equiv v'_i$ mod 2
and thus 
$\,{}^t\!v'(-D\,{}^t\!C)v'\equiv (D\,{}^t\!C)_0v'$ mod $2$, similarly 
$v''(B\,{}^t\!A)v''\equiv (B\,{}^t\!A)_0v''$ mod $2$. 
Next $B\,{}^t\!C\equiv I+A\,{}^t\!D$ mod $2$ and thus $\,{}^t\!v'(D\,{}^t\!A)v''+\,{}^t\!v''(B\,{}^t\!C)v'\equiv \,{}^t\!v'v''$ mod $2$.
Thus we found that
$$
q_\Delta(M^{-1}v)\,=\,\,{}^t\!v'v''+(a\,{}^t\!D-b\,{}^t\!C+(D\,{}^t\!C)_0)v'+
(-a\,{}^t\!B+b\,{}^t\!A+(B\,{}^t\!A)_0)v''
$$
so $q_\Delta(M^{-1}v)=q_{M\cdot \Delta}(v)$, as desired.

\subsection{Transvections}\label{transvections}
The group $Sp(2g,\FF_2)$ is generated by transvections $t_v$, 
for $v\in V$, which are analogous to reflections in orthogonal groups
(\cite{Jac}, $\S$ 6.9).
They are defined as:
$$
t_v:V\longrightarrow V,\qquad t_v(w):=w+E(w,v)v.
$$
It is straightforward to verify that $t_v\in Sp(2g,\FF_2)$, in fact the same formula works also for $\ZZ$ instead of $\FF_2$ and then defines elements in $Sp(2g,\ZZ)$. As
$gt_vg^{-1}=t_{g(v)}$ for $g\in Sp(2g,\FF_2)$ and $v\in V$, the
non-trivial transvections form a conjugacy class. 

For $v\in V$ one verifies easily $t_v^2=1$, the identity on $V$, so 
$t_v$ is an involution, and
$$
t_vt_w=t_wt_v\quad\mbox{if}\quad E(v,w)=0,\qquad
t_vt_wt_v=t_{v+w} \quad\mbox{if}\quad E(v,w)=1
$$
for $v,w\in V$. In particular $(t_vt_w)^2=1$ if $E(v,w)=0$ and $(t_vt_w)^3=1$ if $E(v,w)=1$.

Using these Coxeter relations it is easy to establish a relation between
$Sp(6,\FF_2)$ and $W(E_7)$. 
$$
 
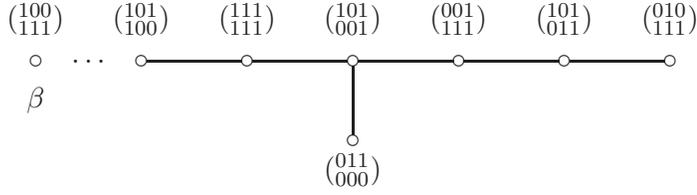
\begin{picture}(250, 95)%
   \put(-40,60){\circle{4}}%
   \put(-20,60){\makebox(0, 0)[c]{$\ldots$}}%
   \put(0, 60){\circle{4}}%
   \put(2, 60){\line(1, 0){36}}%
   \put(40, 60){\circle{4}}%
   \put(42, 60){\line(1, 0){36}}%
   \put(80, 60){\circle{4}}%
   \put(82, 60){\line(1, 0){36}}%
   \put(120, 60){\circle{4}}%
   \put(122, 60){\line(1, 0){36}}%
    \put(160, 60){\circle{4}}%
   \put(162, 60){\line(1, 0){36}}%
   \put(200, 60){\circle{4}}%
   \put(80, 58){\line(0, -1){26}}%
   \put(80, 30){\circle{4}}%
   \put(-40, 70){\makebox(0, 0)[b]{$({}^{100}_{111})$}}%
   \put(0, 70){\makebox(0, 0)[b]{$({}^{101}_{100})$}}%
   \put(40, 70){\makebox(0, 0)[b]{$({}^{111}_{111})$}}%
   \put(80, 70){\makebox(0, 0)[b]{$({}^{101}_{001})$}}%
   \put(120, 70){\makebox(0, 0)[b]{$({}^{001}_{111})$}}%
   \put(160, 70){\makebox(0, 0)[b]{$({}^{101}_{011})$}}%
   \put(200, 70){\makebox(0, 0)[b]{$({}^{010}_{111})$}}%
    \put(80, 25){\makebox(0, 0)[t]{$({}^{011}_{000})$}}%
   \put(-40, 50){\makebox(0, 0)[t]{$\beta$}}%
\end{picture}
$$
One verifies that $E(v,w)=0$ for points $v,w$ in the diagram except when they are connected by an edge (or the dotted edge), in which case $E(v,w)=1$.
Thus there is a surjective homomorphism from $W(E_7)$ to $Sp(6,\FF_2)$.
Note that $\beta$ indicates the longest root of $E_7$ and from the diagram one can see that $\beta^\perp$ contains a root system of type $D_6$. 
This root system can be realized as the set of vectors $\pm(e_i\pm e_j)$ in $\RR^6$, in particular, 
there are 6 mutually perpendicular vectors in $D_6$ (for example $e_1+e_2,e_1-e_2,e_3+e_4,\ldots,e_5-e_6$). 
Therefore there are $7$ perpendicular roots $\beta=\beta_1,\ldots,\beta_7$ in $E_7$.
The product of the reflections in $W(E_7)$ defined by these $7$ roots is obviously $-1$ on $\RR^7$, thus it lies in the center of $W(E_7)$. 
The group $W(E_7)/<-1>$ is simple and is isomorphic to $Sp(6,\FF_2)$,
one actually has
$$
W(E_7)\,\cong\,Sp(6,\FF_2)\,\times\,\{\pm 1\}.
$$
That the image of $-1$ is trivial can also be checked as follows:
the seven perpendicular roots map to the $7$ non-zero points in a
Lagrangian subspace. In suitable coordinates we may assume that this subspace is $L:=\{({}^{abc}_{000}):a,b,c\in\FF_2\}$ and then one verifies easily that the product of the $7$ transvections $t_v$, where $v$ runs over the non-zero elements of $L$, is the identity on $V$.

\subsection{The action of the transvections on the characteristics}
\label{transvections char}
Let $v\in V$ and let $q$ be a quadratic form with associated bilinear
form $E$ and characteristic $\Delta_q$.
As $t_v$ is an involution, $q(v+w)=q(v)+q(w)+E(v,w)$ for all $v,w\in V$ and $q(av)=aq(v)$ for $a\in\FF_2$ we have
$$
(t_v\cdot q)(w)=q(t_v(w))=q(w+E(v,w)v)=q(w)+E(v,w)q(v)+E(v,w)^2.
$$
Hence we get the simple rule:
$$
(t_v\cdot q)(-)\,=\,\left\{
\begin{array}{rcl}
q(-)&\mbox{if}&q(v)=1,\\
q(-)+E(v,-)&\mbox{if}&q(v)=0,
\end{array}
\right.\qquad\mbox{so}\quad
t_{({}^{v'}_{v''})}\cdot q_{[{}^a_b]}\,=\,
q_{[{}^{a+v''}_{b+v'}]}
$$
in case $q_{[{}^a_b]}(({}^{v'}_{v''}))=0$.

Note that the $7$ points in the top row of the Dynkin diagram 
in \ref{transvections}
are a diagram 
of type $A_7$ and that $W(A_7)=S_8$, so we do get a copy of 
$S_8\subset W(E_7)$, which maps, isomorphically, to a subgroup of $Sp(6,\FF_2)$. For each of these seven points 
$v=({}^{abc}_{def})\in V$ one verifies that $q(v):=ad+be+cf=1$
where $q$ is the quadric with characteristic 
$\Delta_q=[{}^{000}_{000}]$. Hence we get seven transvections
in the orthogonal group $O^+$ of $q$ which generate a subgroup isomorphic to $S_8$. As $O^+$ has index $36$ in
$Sp(6,\FF_2)$ it follows that $O^+\cong S_8$.

\section{B: The Heisenberg group and the theta constants}\label{B}

\subsection{The projective representation of $\Gamma_g$ on $T_g$}
\label{proj repr}
Recall that we defined 
$\Theta[\sigma](\tau):=\theta[{}^\sigma_{0}](2\tau)$. Let $T_g$ be
the $2^g$-dimensional vector space spanned by the holomorphic
functions $\Theta[\sigma]$, with $\sigma\in \FF_2^g$:
$$
T_g\,:=\,
\langle\ldots,\Theta[\sigma],\ldots\rangle_{\sigma\in\FF_2^g}.
$$
There is a (projective, linear) action of $\Gamma_g$ on $T_g$
given by
$$
(M\cdot \Theta[\sigma])(\tau)\,:=\,
(\det(C\tau+D))^{-1/2}\Theta[\sigma](M^{-1}\cdot \tau),
$$
the sign ambiguity leads to a projective, rather then a linear, 
representation which factors over the finite group $\Gamma_g/\Gamma_g(2,4)$.
This action obviously induces 
usual action on the $\theta[\Delta]^4\in Sym^4(T_g)$
as well as the $Sp(2g)$-representations $\rho_k$ on the Heisenberg invariants of degree $2k$ in $\CC[\ldots,\Theta[\sigma],\ldots]$.

To compute the action of $M\in\Gamma_g$ on $T_g$ one first considers the
case that $M$ has blocks $A,\ldots,D$ with $C=0$, in which case it is straightforward to find $M\cdot \Theta[\sigma]$.
Next one considers the matrix $S\in\Gamma_g$ given by
$A=D=0$, $B=-C=I$ on the $\Theta[\sigma]$. 
$$
\Theta[\sigma](S^{-1}\cdot \tau)\,=
\,\Theta[\sigma](-\tau^{-1})\,=\,
\theta[{}^\sigma_0](-2\tau^{-1})\,=\,
c_\tau\theta[{}^0_\sigma](\tau/2),
$$
where in the last equality we use the classical transformation formula (cf.\ \cite{Igusa}, II.5, Corollary), in particular $c_\tau$ does not depend on $\sigma$. 
To express $\theta[{}^0_\sigma](\tau/2)$ in terms of the 
$\Theta[\sigma]$'s we write each $m\in\ZZ^g$ uniquely as $m:=2n+\rho$
with $\rho=(\rho_1,\ldots,\rho_g)$, $\rho_j\in\{0,1\}$. Then the sum over $m\in\ZZ^g$ is the sum over $n\in\ZZ^g$ and over these $2^g$
elements $\rho$:
$$
\theta[{}^0_\sigma](\tau/2)\,=\,\sum_\rho\sum_{n\in\ZZ^g}
e^{\pi i((2n+\rho)(\tau/2)\,{}^t\!(2n+\rho)+(2n+\rho)\,{}^t\!\sigma)}
\,=\,\sum_\rho\sum_{n\in\ZZ^g}
e^{\pi i(n+\rho/2)2\tau\,{}^t\!(n+\rho/2)+
2\pi in\,{}^t\!\sigma+\pi i\rho\,{}^t\!\sigma},
$$
as $n\,{}^t\!\sigma\in\ZZ$ we get:
$$
\Theta[\sigma](S^{-1}\cdot\tau)
\,=\,
c_\tau\sum_\rho(-1)^{\rho\sigma}\Theta[\rho](\tau).
$$
This transformation on the $\Theta[\sigma]$'s is also known as the (finite) Fourier transform.

As $\Gamma_g$ is generated by the matrices with $C=0$ and $S$ (cf.\
\cite{Igusa}, I.10, Lemma 15) this allows one in principle to compute
the action of any $M\in \Gamma_g$ on $T_g$.

\subsection{The actions of $H_g$ and $\Gamma_g(2)/\Gamma_g(2,4)$}
\label{actHG}
We show that the Heisenberg group $H_g$ (cf.\ section \ref{Heisenberg})
and
the subgroup $\Gamma_g(2)/\Gamma_g(2,4)$ of $\Gamma_g/\Gamma_g(2,4)$
act in the same way
on the vector space $T_g$ spanned by the $\Theta[\sigma]$'s.
In particular, the Heisenberg invariants in $\CC[\ldots,\Theta[\sigma],\ldots]$ are the $\Gamma_g(2)/\Gamma_g(2,4)$-invariants and these are thus modular forms
on $\Gamma_g(2)$.

First of all we consider the exact sequence:
{\renewcommand{\arraystretch}{1.7}
$$
\begin{array}{ccccccccc}
0&\longrightarrow&\Gamma_g(2)/\Gamma_g(2,4)&\longrightarrow&
\Gamma_g/\Gamma_g(2,4)&\longrightarrow&\Gamma_g/\Gamma_g(2)&
\longrightarrow&0\\
&&\cong\downarrow\phi
&&&&\cong\downarrow\phantom{\cong}&&\\
&&H_g/\mu_4 &&&&Sp(2g,\FF_2)&&
\end{array}
$$
}
where, for $M=I+2S\in \Gamma_g(2)$ and $S$ with blocks $A',B',C',D'$, the isomorphism $\phi$ is induced by
$$
\tilde{\phi}\,:\,\Gamma_g(2)\,\longrightarrow\,\FF_2^g\times\FF_2^g,
\qquad
M=I+2S\,\longmapsto \, ({\rm diag }C',{\rm diag}B')\;{\rm mod}\,2.
$$
It is easy to verify that $\tilde{\phi}$ is a homomorphism, its kernel
consists of the $M$ with 
${\rm diag}B'\equiv {\rm diag }C'\equiv \,0\;{\rm mod}\,2$,
or equivalently, 
${\rm diag}(I+2A')2\,{}^t\!B'
\equiv {\rm diag }2C'\,{}^t\!(1+2D')\equiv \,0\;{\rm mod}\,4$,
so $\ker\tilde{\phi}=\Gamma_g(2,4)$.
Using $M=I+2S$ with $A'=D'=0$ and with either $B'$ diagonal and $C'=0$ or $B'=0$ and $C'$ diagonal (such matrices are symplectic) one verifies the surjectivity of $\tilde{\phi}$.

Let $M\in\Gamma_g(2)$ be the matrix with blocks $A=D=I$, $B={\rm diag}(2u_1,\ldots,2u_g)$ with $u_i\in\{0,1\}$ and $C=0$.
Then one easily verifies:
$$
(M\cdot \Theta[\sigma])(\tau)\,=\,
\Theta[\sigma](M^{-1}\cdot \tau)\,=\,\Theta[\sigma](\tau-B)\,=\,
(-1)^{\sigma u}\Theta[\sigma](\tau),\qquad u:=(u_1,\ldots,u_g).
$$
In particular, $M$ acts as $(1,0,u)\in H_g$ on $T_g$.
Using the identity
$$
\left(\begin{array}{cc} I&0\\C&I\end{array}\right)\;=\;
\left(\begin{array}{cc} 0&I\\-I&0\end{array}\right)
\left(\begin{array}{cc} I&-C\\0&I\end{array}\right)
\left(\begin{array}{cc} 0&-I\\I&0\end{array}\right),
$$
one checks
that if $N\in\Gamma_g(2)$ is the matrix with blocks $A=D=I$, $B=0$, $C=2{\rm diag}(x_1,\ldots,x_g)$,
then
$N\cdot \Theta[\sigma]$ is $\Theta[\sigma+x]$, up to a multiplicative factor independent of $\sigma$. 
In particular, $N$ acts as $(1,x,0)\in H_g$ on $T_g$.

Summarizing, if $M\in \Gamma_g(2)$ and $\phi(M)=(x,u)\in\FF_2^{2g}$,
then the action of $M$ and $(t,x,u)\in H_g$ (for any $t\in\mu_4$) 
on $\PP T_g$ coincide.
Thus the Heisenberg group action on $T_g$ is a lifting 
of the action of $\Gamma_g(2)/\Gamma_g(2,4)$ on $\PP T_g$.

\subsection{The action of a transvection on the $\Theta[\sigma]$'s}
\label{transvections theta}
For $v\in \FF_2^{2g}=\ZZ^{2g}/2\ZZ^{2g}$, $v\neq 0$, we defined the transvection
$t_v\in Sp(2g)$ in \ref{transvections}.
Let $\tilde{v}\in\ZZ^{2g}$ a vector such that $\tilde{v}$ maps to
$v\in \ZZ^{2g}/2\ZZ^{2g}$ and let $t_{\tilde{v}}\in \Gamma_g$ be the associated transvection:
$$
t_{\tilde{v}}\,:\,\ZZ^{2g}\longrightarrow \ZZ^{2g},\qquad
x\longmapsto x+E(x,\tilde{v})\tilde{v},
$$ 
where now $E$ is the standard symplectic form on $\ZZ^{2g}$ as in 
\ref{sp char}. Then the image of $t_{\tilde{v}}$ in $Sp(2g)=\Gamma_g/\Gamma_g(2)$ is $t_v$. In particular, the linear maps
$\rho_k(t_v)$ on the modular forms of weight $2k$ on $\Gamma_g(2)$
are induced by the action of $t_{\tilde{v}}$
on $T_g$ as in \ref{proj repr}:
$$
t_{\tilde{v}}\cdot f\,=\,\rho_k(t_v)f\qquad 
f\in \CC[\ldots,\Theta[\sigma],\ldots]_{4k}^{H_g}
=M_{2k}^\theta(\Gamma_g(2)).
$$
We give a simple description of the action of $t_{\tilde{v}}$. 

Let $v=(x,u)$ and let, as in 
section \ref{Heisenberg}, $xu=x_1u_1+\ldots+x_gu_g\in\FF_2$. 
Define $U_v:T_g\rightarrow T_g$ to be the action of $(1,v)\in H_g$ 
in case $xu=1$ and let it denote the action of $(i,v)$ (with $i^2=-1$)
in case $xu=0$. In particular, we have $U_v^2=-I$ for all $v\in V$.
Next we define a linear map on $T_g$ by:
$$
\tilde{t}_v\,:\,T_g\,\longrightarrow\,T_g,\qquad
\tilde{t}_v\,:=\,\frac{1-i}{2}(U_v+I).
$$
We will now verify that $\tilde{t}_v=t_{\tilde{v}}$ on $\PP T_g$
and that $\tilde{t}_v$ induces $\rho(t_v)$ on the Heisenberg invariants.

\subsubsection{Conjugation by $\tilde{t}_v$.}\label{conjugation}
An important property of $\tilde{t}_v$ is that for all $w\in V$ there is a non-zero constant $c_{v,w}\in\CC$
such that:
$$
\tilde{t}_vU_w\tilde{t}_v^{-1}\,=\,c_{v,w}U_{t_v(w)}.
$$
For example, if $xu=0$, then $U_v^{-1}=-U_v$ and so
$\tilde{t}_v^{-1}=\frac{1+i}{2}(-U_v+I)$, hence
$$
\tilde{t}_vU_w\tilde{t}_v^{-1}\,=\,
\mbox{$\frac{1}{2}$}(U_vU_wU_v^{-1}+U_vU_w-U_wU_v+U_w),
$$
in case $E(v,w)=0$ one has $U_vU_w=U_wU_v$ hence $\tilde{t}_vU_w\tilde{t}_v^{-1}=U_w$ and if $E(v,w)=1$ one has $U_vU_w=-U_wU_v$
which is also equal to $U_{v+w}$ up to a constant,
hence $\tilde{t}_vU_w\tilde{t}_v^{-1}=(U_vU_w-U_wU_v)/2=cU_{v+w}$.
Similarly one handles the case $xu=1$.

\subsubsection{The action of the lifts on $T_g$.}
As $\Gamma_g(2)/\Gamma_g(2,4)$ is a normal subgroup
of $\Gamma_g(1)/\Gamma_g(2,4)$, the group $\Gamma_g(1)/\Gamma_g(2,4)$ acts by conjugation on $\Gamma_g(2)/\Gamma_g(2,4)$. 
In particular one verifies easily that for $\tilde{v},\tilde{w}\in\ZZ^{2g}$ we have
$$
t_{\tilde{v}}\,t_{2\tilde{w}}\,t_{\tilde{v}}^{-1}\,=\,
t_{2\tilde{w}'},\qquad \mbox{with}\quad\tilde{w}'=t_{\tilde{v}}(\tilde{w}).
$$
In \ref{actHG} we verified that $t_{2\tilde{w}}\in \Gamma_g(2)/\Gamma_g(2,4)$ acts as $U_w$
on $\PP T_g$. Thus
$t_{\tilde{v}}\in Sp(2g,\ZZ)$ must act on $\PP(T_g)$
by a linear map $S_{\tilde{v}}$ which satisfies 
$S_{\tilde{v}}U_wS_{\tilde{v}}^{-1}=U_{t_v(w)}$.

On the other hand, we verified in \ref{conjugation} that $\tilde{t}_v$ has this property. 
Thus $S:=\tilde{t}_v^{-1}S_{\tilde{v}}$ satisfies $SU_v=U_vS$ for all $v\in V$. But the group generated by the $U_v$ acts irreducibly on
$T_g$ (this is easy to check). Therefore, by Schur's Lemma,
$S$ is a scalar multiple of the identity and thus 
$$
\tilde{t}_v (\Theta[\sigma])\,=\,t_{\tilde{v}}\cdot \Theta[\sigma],
\qquad(\forall \sigma)
$$
up to a scalar independent of $\sigma$.

\subsubsection{The action of $Sp(2g)$ on the Heisenberg invariants.}
\label{action on invariants}
Now we show that the action of $t_v$ on the Heisenberg invariants
in $\CC[\ldots,\Theta[\sigma],\ldots]$ is given by $\tilde{t}_v$. 
We already know that the action of $\tilde{t}_v$
is the action of $t_v$ up to a scalar multiple. 
As the degree of a Heisenberg invariant polynomial is a multiple of four, the action of $\tilde{t}_v$
and $i^a \tilde{t}_v$, for integer $a$, are the same on the Heisenberg invariants.

First of all we note that $\tilde{t}_v^2=iU_x$ and thus $\tilde{t}_v^2$ acts
as the identity on the Heisenberg invariants. Also $t_v$ has order two, so if
$\lambda_v \tilde{t}_v$ gives the action of $t_v$ then $\lambda_v^8=1$.
To exclude the case that $\lambda_v^4=-1$, we consider the action of
$\tilde{t}_v$ on the $\theta[\Delta]^4$. This action should coincide
with the one of $t_{\tilde{v}}\in \Gamma_g$. 
Choose an even characteristic $\Delta$ such that 
$t_{\tilde{v}}\cdot\Delta=\Delta$. 
As the $t_v$ generate the orthogonal group of the corresponding quadratic form $q_\Delta$, we have $\rho_2(t_{\tilde{v}})\theta[\Delta]^4=-\theta[\Delta]^4$.

Now an easy computation shows that, with
the definition of $\tilde{t}_v$ above
and the classical formula for $\theta[\Delta]^2$ given in \ref{classical},
one has:
$$
\tilde{t}_v(\theta[\Delta]^2)\,=\,\pm i\theta[\Delta]^2\qquad
\mbox{if}\quad t_{\tilde{v}}\cdot\Delta=\Delta.
$$
Hence also $\tilde{t}_v(\theta[\Delta]^4)=-\theta[\Delta]^4$.
Therefore $\tilde{t}_v$ gives the action of $t_v$ on the Heisenberg invariants.

\subsubsection{Example.} 
Let $g=1$ and let $v=({}^0_1)\in \FF_2^2$. Then w.r.t.\ to the basis $\Theta[0],\Theta[1]$ of $T_1$ we have
$$
U_v\,=\,i\begin{pmatrix} 1&0\\0&-1 \end{pmatrix},\qquad
\tilde{t}_v\,=\,\frac{1-i}{2}(U_v+I)\,=\,
\begin{pmatrix} 1&0\\0&-i \end{pmatrix}\qquad(\in GL(T_1)).
$$
From this we find the action for example on $\theta[{}^0_0]$:
$$
\tilde{t}_v(\theta[{}^0_0]^4)\,=\,
\tilde{t}_v(\Theta[0]^2+\Theta[1]^2)^2\,=\,
(\Theta[0]^2-\Theta[1]^2)^2=\theta[{}^0_1]^4.
$$
Let $\tilde{v}=({}^0_1)\in \ZZ^2$, then 
$$
t_{\tilde{v}}\,=\,t_{({}^0_1)}\,=\,
\begin{pmatrix} 1&1\\0&1 \end{pmatrix}=T\qquad(\in \Gamma_1)
$$
and thus we recover the action of $T$ on $\theta[{}^0_0]^4$.

Note that if $v=(0,u)\in\FF_2^g\times\FF_2^g$ then
$U_v$ is diagonal, and so is $\tilde{t}_v$ (which corresponds to an element with $C=0, A=D=I$ in $Sp(2g,\ZZ)$), so we can compute the action in a particularly simple (and fast!) manner.

\end{appendix}

\bibliographystyle{my-h-elsevier}

\end{document}